\documentclass{article}

\usepackage{a4wide, amsthm, amssymb, amscd}
\usepackage[all]{xy}
\usepackage{longtable}

\newtheorem{theorem}{Theorem}[section]
\newtheorem{lemma}[theorem]{Lemma}
\newtheorem{proposition}[theorem]{Proposition}
\newtheorem{corollary}[theorem]{Corollary}
\newtheorem{remark}[theorem]{Remark}
\newtheorem{definition}[theorem]{Definition}

\newtheorem{example}[theorem]{Example}

\newtheorem{assumption}[theorem]{Assumption}
\newtheorem{notation}[theorem]{Notation}

\newcounter{namedenum}
\newenvironment{namedenumerate}[1][]{%
  \begin{list}{\textup{(#1\arabic{namedenum})}}%
    {\settowidth{\leftmargin}{\textup{(#1\arabic{namedenum})}}
     \addtolength{\leftmargin}{1em}
     \settowidth{\labelwidth}{\textup{(#1\arabic{namedenum})}}
     \addtolength{\labelwidth}{1em}
     \setlength{\itemindent}{0em}
     \setlength{\labelsep}{0.6em}
     \usecounter{namedenum}}}%
  {\end{list}}

\newcommand{\p}{\mathfrak p}
\newcommand{\rp}{\p^\Lsh}

\newcommand{\s}{\mathfrak s}
\newcommand{\rs}{{\s^\Lsh}}

\newcommand{\rwedge}{\wedge^\Lsh}

\newcommand{\rvee}{\vee^\Lsh}
\newcommand{\rt}[1]{{#1}^{\Lsh}}
\newcommand{\riota}{\iota^{\Lsh}}
\newcommand{\rphi}{\varphi^{\Lsh}}
\newcommand{\NN}{\mathbb N}
\newcommand{\ZZ}{\mathbb Z}

\newcommand{\wordlen}[1]{\parallel\!\!#1\!\!\parallel}
\newcommand{\sinf}{{\rm\inf_s}}
\newcommand{\ssup}{{\rm\sup_s}}

\newcommand{\CX}{x^G}

\newcommand{\SSS}{\textup{SSS}}

\newcommand{\USS}{\textup{USS}}
\newcommand{\SC}{\textup{SC}}
\newcommand{\SCG}{\textup{SCG}}

\title{Solving the conjugacy problem in Garside groups \\ by cyclic sliding}
\author{
  Volker Gebhardt%
      \footnote{Both authors partially supported by MTM2007-66929 and FEDER.}
\and
  Juan Gonz\'{a}lez-Meneses$^{\ast,}$%
      \footnote{This work was done partially while the second author was visiting the Institute for
      Mathematical Sciences, National University of Singapore in 2007. The visit was supported by
      the Institute.}
}
\date{September 4, 2008}

\begin{document}

\maketitle


\addtolength{\topmargin}{-15mm}

\begin{abstract}
We present a solution to the conjugacy decision problem and the conjugacy search problem in Garside
groups, which is theoretically simpler than the usual one, with no loss of efficiency. This is done
by replacing the well known cycling and decycling operations by a new one, called cyclic sliding,
which appears to be a more natural choice.

We give an analysis of the complexity of our algorithm in terms of fundamental operations with
simple elements, so our analysis is valid for every Garside group. This paper intends to be
self-contained, not requiring any previous knowledge of prior algorithms, and includes all the
details for the algorithm to be implemented on a computer.
\end{abstract}

\section{Introduction}

The {\bf Conjugacy Decision Problem} (CDP) for a group $G$ is the decision problem of determining, given any two elements $a,b\in G$, whether $a$ and $b$ are conjugate in $G$.
The {\bf Conjugacy Search Problem} (CSP), on the other hand, requires to compute for any two given conjugate elements $a,b\in G$ a conjugating element $c$ such that $c^{-1} ac=b$. (We will also write $a^c=b$.)

In this paper we will describe a new algorithm to solve both problems in Garside groups
(of finite type). The simplicity of the algorithm will allow us to describe it completely in this
introduction in a ready-to-implement manner. The main difference to established algorithms is the
use of an operation called {\em cyclic sliding}, which is a special kind of conjugation introduced in~\cite{GG1}.  Cyclic sliding assumes the role played by {\it cycling} and {\it decycling} in previous algorithms.

Cyclic sliding will be motivated and explained in \S\ref{SS:new ingredients}, but it can be defined
right now. One just needs to recall the following notions in a Garside group $G$, which are well
known to specialists. Firstly, $G$ admits a partial order $\preccurlyeq$, and there is a special
element $\Delta$, called Garside element. Given  $x\in G$, $\inf(x)$ and $\sup(x)$ are the maximal
and minimal integers, respectively, satisfying
$\Delta^{\inf(x)}\preccurlyeq x \preccurlyeq\Delta^{\sup(x)}$.
Secondly, given $a,b\in G$, there is a unique greatest common divisor $a\wedge b$ with respect to
$\preccurlyeq$. Finally, the elements in the set
$[1,\Delta]=\{s\in G\ |\ 1\preccurlyeq s\preccurlyeq \Delta\}$, called {\it simple elements}, generate $G$.
We assume this set to be finite (that is, $G$ is of finite type). It is well known how to
compute all the above data in a Garside group $G$ of finite type, as we shall see.

\pagebreak

Using the above well known notions, we can define the following:

\begin{definition}[\cite{GG1}]
Given $x\in G$, we define the {\bf preferred prefix} $\p(x)$ of $x$ as the
simple element
$$
\p(x)=\left(x\Delta^{-\inf(x)}\right) \wedge
\left(x^{-1}\Delta^{\sup(x)}\right)\wedge \Delta,
$$
and we define the {\bf cyclic sliding} $\s(x)$ of $x$ as the conjugate
of $x$ by its preferred prefix, that is,
$$
\s(x)= x^{\p(x)}= \p(x)^{-1} x \: \p(x).
$$
\end{definition}

This is enough to describe a simple algorithm to solve the conjugacy decision problem and the
conjugacy search problem in a Garside group of finite type. The algorithm we present now, however,
is by far not the best possible one. In \S\ref{SS:algorithm} we will give a much better
algorithm, which requires some other notions besides the preferred prefix and the cyclic sliding.
Nevertheless, the simple version given here for illustration can be useful for theoretical purposes
or for applying it to small examples.

\bigskip

\begin{minipage}{\textwidth}
\fbox{{\sc Algorithm 0:} {\bf Solving the conjugacy problem in a Garside group $G$ of finite type}}

\begin{longtable}{lp{12.5cm}}
{\bf Input:} & $x,y\in G$. \\
{\bf Output:} & - Whether $x$ and $y$ are conjugate.\newline
                - If $x$ and $y$ are conjugate, an element $c$ such
                      that $x^c=y$.
\end{longtable}
\begin{enumerate}
\item Set $\widetilde x=x$, $c_1=1$ and $\mathcal T =\emptyset$.
\item While $\widetilde x\notin \mathcal T$, set
  $\mathcal T=\mathcal T \cup \{\widetilde x\}$,
  \ $c_1=c_1 \cdot \p(\widetilde x)$ \ and
  \ $\widetilde x=\s(\widetilde x)$.
\item Set $\widetilde y=y$, $c_2=1$ and $\mathcal T =\emptyset$.
\item While $\widetilde y\notin \mathcal T$, set $\mathcal
  T=\mathcal T \cup \{\widetilde y\}$,\ $c_2=c_2 \cdot
  \p(\widetilde y)$ \ and  \ $\widetilde y=\s(\widetilde y)$.
\item  Set $\mathcal V=\{\widetilde x\}$, $\mathcal V'=\{\widetilde x\}$ and
  $c_{\widetilde x}=1$.
\item While $\mathcal V' \neq \emptyset$, do:
   \begin{enumerate}
   \item Take $v\in \mathcal V'$.
   \item For every simple element $s$, do:
      \begin{enumerate}
      \item If $v^s=\widetilde y$, then set $c_{\widetilde y}=c_v\cdot s$.
            Return `$x$ and $y$ are conjugate by $c_1 \cdot c_{\widetilde y} \cdot c_2^{-1}\;$'.
            STOP.
      \item If $v^s\notin \mathcal V$, then:
          \begin{enumerate}
          \item  Apply iterated cyclic sliding to $v^s$ until the first repetition, say $w$.
          \item  If $w=v^s$, then set $c_{v^s}=c_v\cdot  s$, $\mathcal V=\mathcal V\cup \{v^s\}$,
                 and $\mathcal V'=\mathcal V'\cup \{v^s\}$.
          \end{enumerate}
      \end{enumerate}
   \item Remove $v$ from $\mathcal V'$.
   \end{enumerate}
\item Return `$x$ and $y$ are not conjugate'.
\end{enumerate}
\end{minipage}
\bigskip


The set $\mathcal V$ computed by the above algorithm, called the {\it set of sliding circuits of
$x$} and denoted $\SC(x)$, was introduced in~\cite{GG1}. It is a finite invariant of the conjugacy
class $\CX$ of $x$, that is, it is a finite subset of $\CX$ and only depends on $\CX$, not on $x$
itself. This set $\SC(x)$ consists of those conjugates of $x$ which are stabilised by $\s^k$ for
some positive integer $k$ and it is analogous to the ultra summit set $\USS(x)$
from~\cite{Gebhardt}. One has $\SC(x)\subseteq \USS(x)$, and in general $\SC(x)$ is a proper subset
of $\USS(x)$.

The first two lines of the algorithm compute an element $\widetilde x \in \SC(x)$, by applying
iterated cyclic sliding until the first repetition is reached (which is $\widetilde x$). A
conjugating element $c_1$ from $x$ to $\widetilde x$ is also computed. The following two lines
compute $\widetilde y\in \SC(y)$ and a conjugating element $c_2$ from $y$ to $\widetilde y$ in
the same way. Then, the algorithm starts to compute the whole set $\SC(x)$. If during the
computation it finds $\widetilde y$ as an element of $\SC(x)$, the algorithm stops and returns a
conjugating element from $x$ to $y$. If this does not occur, that is, if the algorithm computes the
whole set $\SC(x)$ without finding $\widetilde y$ in it, then it returns the message `$x$ and $y$
are not conjugate'.

The use of cyclic sliding not only allows to develop a simpler algorithmic solution to the CDP/CSP,
but also is of theoretical interest; we refer to~\cite{GG1} for details. It is shown there that the
set of sliding circuits has all the good properties of the ultra summit set, but is the more
natural invariant in many ways.  In particular, the properties of the set of sliding circuits fully
extend to the case of elements of {\it summit canonical length} 1, which is not the case for ultra
summit sets.  Another indication of the naturalness of the cyclic sliding operation is the fact
that for {\it super summit elements} which have a {\it rigid} conjugate, the (unique) minimal
positive element yielding a rigid conjugate is precisely the conjugating element obtained by
iterated cyclic sliding.

The structure of this paper is as follows. In the introduction, we present our algorithm solving
the conjugacy problems in Garside groups in a ready-to-implement form.  This presentation is kept
as concise as possible; explanations, motivations and the proof of correctness are postponed to
later sections. More precisely, in \S\ref{SS:Garside}, we give a basic introduction to the theory
of Garside groups; specialists may skip this part. In \S\ref{SS:new ingredients} we briefly explain
the new concepts from~\cite{GG1} which are subsequently used for the detailed description of the
algorithm in \S\ref{SS:algorithm}.

The rest of the paper is devoted to the explanation and analysis of the algorithm.
\S\ref{S:cyclic sliding} contains a summary of results from~\cite{GG1} which are required in our discussion.  In \S\ref{S:description of the algorithm} the algorithm is explained and shown to be
correct.  Finally, the complexity of the new algorithm is analysed in \S\ref{S:complexity}, where
\S\ref{S:ComputingInGarsideGroups} discusses how the operations required for our algorithm can be
realised, only assuming knowledge of the lattice of simple elements.

\subsection{Basic facts about Garside groups}\label{SS:Garside}

Garside groups were defined by Dehornoy and Paris~\cite{DP}. For a detailed
introduction to these groups, see~\cite{Dehornoy}; a shorter introduction,
containing all the details needed for this paper can be found
in~\cite{BGG1} (\S 1.1 and the beginning of \S 1.2).

One of the possible definitions of a Garside group is the following. A
group $G$ is said to be a  {\bf Garside group} with
{\bf Garside structure $(G,P,\Delta)$} if it admits a submonoid $P$
satisfying $P\cap P^{-1}=\{1\}$, called the monoid of
{\bf positive elements}, and a special element $\Delta\in P$ called
the {\bf Garside element}, such that the following properties hold:
\begin{namedenumerate}[G]\vspace{-\topsep}
\item The partial order $\preccurlyeq$ defined on $G$ by $a\preccurlyeq
  b \Leftrightarrow a^{-1}b\in P$ (which is invariant under
  left multiplication by definition) is a lattice order. That is, for every
  $a,b\in G$ there are a unique least common multiple $a\vee b$ and a
  unique greatest common divisor $a\wedge b$ with respect to $\preccurlyeq$.
\item The set $[1,\Delta]=\{a\in G\ |\ 1\preccurlyeq a \preccurlyeq
  \Delta\}$, called the set of {\bf simple elements}, generates $G$.
\item Conjugation by $\Delta $ preserves $P$ (so it preserves the
  lattice order $\preccurlyeq$). That is, $\Delta^{-1}P\Delta =P$.
\item For all $x\in P\backslash\{1\}$, one has:
$$
||x|| = \sup\{k \ | \ \exists\, a_1,\ldots,a_k \in P\backslash \{1\} \mbox{ such that }
  x=a_1\cdots a_k \} < \infty.
$$
\end{namedenumerate}

\medskip

\begin{definition}
A Garside structure $(G,P,\Delta)$ is said to be {\bf of finite type} if the
set of simple elements $[1,\Delta]$ is finite.
A group $G$ is said to be a {\bf Garside group of finite type} if it admits a
Garside structure of finite type.
\end{definition}

Throughout this paper, let $G$ be a Garside group of finite type with a fixed Garside structure
$(G,P,\Delta)$ of finite type.

{\bf Remarks:}
\begin{enumerate}\vspace{-\topsep}
\item By definition, $p\in P \Leftrightarrow 1\preccurlyeq p$.
  Given two positive elements $a\preccurlyeq b$, one usually says that $a$ is a \textbf{prefix}
  of $b$.
  Hence the simple elements are the positive prefixes of $\Delta$.

\item The number $||x||$ defined above for each $x\in P\backslash \{1\}$,
  defines a norm in $P$ (setting $||1||=0$). Note that the existence of this
  norm implies that every element in $P\backslash\{1\}$ can be written as a
  product of {\bf atoms}, where an atom is an element $a\in P$ that cannot be
  decomposed in $P$, that is, $a=bc$ with $b,c\in P$ implies that either
  $b=1$ or $c=1$. In any decomposition of $x$ as a product of $||x||$ factors
  in $P\backslash\{1\}$, all of them are atoms. Notice that the set of atoms
  generates $G$. Moreover, the set of atoms is finite if $G$ is of finite type.
\end{enumerate}

The main examples of Garside groups of finite type are Artin-Tits groups of
spherical type. In particular, braid groups are Garside groups. In the braid
group $B_n$ on $n$ strands with the  usual Garside structure that we call
{\bf Artin Garside structure} of $B_n$, one has the following:
\begin{itemize}\vspace{-\topsep}

  \item The atoms are the standard generators $\sigma_1,\ldots,\sigma_{n-1}$.

  \item The positive elements are the braids that can be written as a word
    which only contains positive powers of the atoms.

  \item The simple elements are the positive braids in which any two strands
    cross at most once. One has $|[1,\Delta]|=n!$, so this is a finite type
    Garside structure.

  \item The Garside element $\Delta$ is the positive braid in which any two strands cross exactly
    once (also called the {\em half twist}). That is,
    $\Delta=\sigma_1 (\sigma_2\sigma_1) (\sigma_3\sigma_2\sigma_1)\cdots
    (\sigma_{n-1}\cdots \sigma_1).$

\end{itemize}

Note that the monoid $P$ induces not
only a partial order~$\preccurlyeq $ which is invariant under left multiplication, but
also a partial order $\succcurlyeq $ which is invariant under right multiplication. The latter
is defined by $a\succcurlyeq b \Leftrightarrow ab^{-1}\in P$.  It is obvious from the definitions
that $a\preccurlyeq b$ is equivalent to $a^{-1}\succcurlyeq b^{-1}$.  It follows from
the properties of $G$ that $\succcurlyeq$ is also a lattice order, that $P$ is
the set of elements $a$ such that $a\succcurlyeq 1$, and that the simple
elements are the positive suffixes of $\Delta$ (where we say that a positive
element $b$ is a suffix of $a$ if $a\succcurlyeq b$). We will denote by
$x\rwedge y$ (resp.~$x\rvee y$) the greatest common divisor (resp.~least
common multiple) of $x,y\in G$ with respect to $\succcurlyeq$.

The following notions are well known to specialists in Garside groups:

\begin{definition}\label{D:complement}
Given a simple element $s$, the {\bf right complement} of $s$ is defined by
$\partial(s)=s^{-1}\Delta$, and the {\bf left complement} of $s$ is
$\partial^{-1}(s)=\Delta\: s^{-1}$.
\end{definition}

Notice that the map $\partial: [1,\Delta] \rightarrow [1,\Delta]$ is a
bijection of the (finite) set $[1,\Delta]$. Notice also that
$\partial^2(s)=\Delta^{-1}s \Delta$. We denote by $\tau$ the inner
automorphism of $G$ corresponding to conjugation by~$\Delta$. Hence
$\partial^2(s)=\tau(s)$.

\begin{definition}\label{D:weighted}
Given two simple elements $a$ and $b$, we say that the decomposition $a\cdot b$ is
{\bf left weighted} if $\partial (a) \wedge b=1$ or, equivalently, if
$ab\wedge \Delta=a$. We say that the decomposition $a\cdot b$ is
{\bf right weighted} if $a \rwedge \partial^{-1}(b)=1$ or, equivalently, if
$ab\rwedge \Delta = b$.

The process of bringing a product $a\cdot b$ of two two simple elements $a$ and $b$ into left
weighted form by replacing it with the product $(as)\cdot (s^{-1}b)$, where
$s=\partial(a)\wedge b$, is called a {\bf local left sliding} or simply a
{\bf local sliding}~\cite{GG1}.  {\bf Local right sliding} is defined analogously.
\end{definition}

\begin{definition}\label{D:normal form}
Given $x\in G$, we say that a decomposition $x=\Delta^p x_1\cdots x_r$, where
$p\in \mathbb Z$ and $r\geq 0$, is the {\bf left normal form} of $x$ if
$x_i\in [1,\Delta]\backslash \{1,\Delta\}$ for $i=1,\ldots,r$ and $x_ix_{i+1}$
is a left weighted decomposition for $i=1,\ldots,r-1$.  We say that a
decomposition $x=y_1\cdots y_r \Delta^p$ is the {\bf right normal form} of $x$
if $y_i\in [1,\Delta]\backslash \{1,\Delta\}$ for $i=1,\ldots,r$ and
$y_iy_{i+1}$ is a right weighted decomposition for $i=1,\ldots,r-1$.
\end{definition}

It is well known that left and right normal forms of elements in $G$ exist and are unique.
(Proposition~\ref{P:normal_form} recalls how to compute them based on local slidings.)
Moreover, the numbers $p$ and $r$ do not depend on the normal form (left or right) that we are 
considering.

\begin{definition}\label{D:inf_sup_ell}
Given $x\in G$, whose left normal form is $\Delta^p x_1\cdots x_r$ and whose right normal form is
$y_1\cdots y_r \Delta^p$, we define the {\bf infimum}, {\bf canonical length} and {\bf supremum} of
$x$, respectively, by $\inf(x)=p$, $\ell(x)=r$ and $\sup(x)=p+r$.
\end{definition}

It is shown in~\cite{EM} that $\inf(x)$ and $\sup(x)$ are precisely the maximal and minimal
integers, respectively, such that $\Delta^{\inf(x)}\preccurlyeq x \preccurlyeq \Delta^{\sup(x)}$
(or, equivalently, $\Delta^{\sup(x)}\succcurlyeq x \succcurlyeq \Delta^{\inf(x)}$).
Moreover, if $x=\Delta^p x_1\cdots x_r$ is in left normal form, then the left normal form of $x^{-1}$ is precisely $x^{-1}=\Delta^{-(p+r)}\;\partial^{-2(p+r)+1}(x_r)
\;\partial^{-2(p+r-1)+1}(x_{r-1})\cdots \partial^{-2(p+1)+1}(x_1)$.
An analogous relation holds for the right normal forms of $x$ and $x^{-1}$.
This implies in particular that $\inf(x^{-1})=-\sup(x)$, $\sup(x^{-1})=-\inf(x)$ and
$\ell(x^{-1})=\ell(x)$.

\begin{definition}
\label{D:initial_final_factor}
Given $x\in G$ , its {\bf (left) initial factor} $\iota(x)$ is defined as
$\iota(x)=x\Delta^{-\inf(x)}\,\wedge\,\Delta$ and its {\bf (left) final factor}
$\varphi(x)$ is defined as $\varphi(x)=(\Delta^{\sup(x)-1}\wedge x)^{-1}\,x$.
\end{definition}

We remark that if $\ell(x)=r>0$, and $\Delta^p x_1\cdots x_r$ is the left
normal form of $x$, then $\iota(x)=\tau^{-p}(x_1)$ and $\varphi(x)=x_r$. This
explains the names given to these simple elements. Notice also that if $r=0$,
that is, if $x=\Delta^p$, then $\iota(x)=1$ and $\varphi(x)=\Delta$.
From the relation between the normal forms of $x$ and $x^{-1}$, we see that
$\iota(x^{-1})=\partial(\varphi(x))$.
Right versions $\riota(x)=\Delta^{-\inf(x)}x\,\rwedge\,\Delta$ and
$\rphi(x) = x\,(\Delta^{\sup(x)-1}\rwedge x)^{-1}$ can be defined analogously.

\begin{definition}\label{D:SSS}
Let $\CX$ denote the conjugacy class of $x$ in $G$ and define
the {\bf summit infimum} $\sinf(x)$ respectively the {\bf summit supremum} $\ssup(x)$ of $x$ by
$\sinf(x) = \max\{\inf(y)\ |\ y\in \CX\}$ and $\ssup(x) = \min\{\sup(y)\ |\ y\in \CX\}$.
The set $\SSS(x) = \{ y\in \CX\ |\ \inf(x)=\sinf(x),\ \sup(y)=\ssup(x)\}$ is called the
{\bf super summit set} of $x$; the elements of $\SSS(x)$ are called {\bf super summit elements}.
\end{definition}

It is well known that $\SSS(x)\subset \CX$ is non-empty and finite \cite{EM} and it is clear from
the definition that $\SSS(x)$ only depends on the conjugacy class of $x$.  Since, by the above remark, $\inf(y^{-1})=-\sup(y)$ and $\sup(y^{-1})=-\inf(y)$ for all $y\in G$, one has $y\in\SSS(x)$
if and only if $y^{-1}\in\SSS(x^{-1})$.

\subsection{Cyclic sliding}\label{SS:new ingredients}

Before explaining our algorithm, we need to describe the underlying operation called
{\it cyclic sliding} introduced in~\cite{GG1}.
The use of cyclic sliding (instead of the well known {\it cycling} and {\it decycling} operations) is what distinguishes the new algorithm from previously known ones.
The cyclic sliding operation will be motivated and explained in more detail in the following
section.
Here we just give the technical definitions, so that they can be used in the algorithm.
Recall that $G$ is a Garside group of finite type with a fixed finite type Garside structure
$(G,P,\Delta)$.

\begin{definition}\label{D:preferred prefix}
Given $x\in G$, the {\bf preferred prefix} $\p(x)$ of $x$ is the simple element
$$
\p(x)=\left(x\Delta^{-\inf(x)}\right)
  \wedge\left(x^{-1}\Delta^{\sup(x)}\right)\wedge \Delta
 = \iota(x)\wedge\iota(x^{-1}) = \iota(x)\wedge\partial(\varphi(x)),
$$
and the {\bf preferred suffix} $\rp(x)$ of $x$ is the simple element
$$
\rp (x)= \left(\Delta^{-\inf(x)}x\right)
  \rwedge \left(\Delta^{\sup(x)}x^{-1}\right)\rwedge \Delta
 = \riota(x)\rwedge\riota(x^{-1}) = \riota(x)\rwedge\partial^{-1}(\rphi(x)).
$$
\end{definition}

\begin{definition}\label{D:cyclic sliding}
Given $x\in G$, the {\bf cyclic left sliding} $\s(x)$ of $x$
is the conjugate of $x$ by its preferred prefix, that is,
$$
\s(x)= x^{\p(x)}= \p(x)^{-1} x \: \p(x),
$$
and the {\bf cyclic right sliding} $\rs(x)$ of $x$ is the
conjugate of $x$ by the inverse of its preferred suffix:
$$
\rs(x)= x^{\rp(x)^{-1}}= \rp(x) \; x \; \rp(x)^{-1}.
$$
If there is no possible confusion, we will call $\s(x)$ the
{\bf cyclic sliding}, or just the {\bf sliding} of $x$.
\end{definition}

It will be convenient to display conjugations in a graph-theoretical
style. In this way, we shall write $u\stackrel{s}{\longrightarrow} v$
if $u^s=v$ for some $u,s,v\in G$. Hence we have:
$$
\begin{CD}
  x  @>\p(x)>>  \s(x)
\end{CD}
\hspace{1cm} \mbox{and} \hspace{1cm}
\begin{CD}
  x  @<\rp(x)<<  \rs(x).
\end{CD}
$$

Elements for which the preferred prefix (or the preferred suffix) is trivial behave particularly nicely in may ways.

\begin{definition}\label{D:rigid}
An element $x\in G$ is called {\bf left rigid} or just {\bf rigid} if $\p(x)=1$.  Similarly, $x$ is
called {\bf right rigid} if $\rp(x)=1$.
\end{definition}

The concept of rigidity was introduced in~\cite{BGG1} and some of the properties of rigid elements were analysed there.  It is obvious from the definition that left (respectively right) rigid elements are fixed points for left (respectively right) cyclic sliding.  The converse clearly is not true.

The main idea of our algorithm is the following: Iterated application of cyclic sliding sends any element $x\in G$ to a finite subset of its conjugacy class $\CX$.
This subset only depends on $\CX$ and is, in general, small.
Hence, it can be used to solve the CDP and the CSP efficiently.
This set is defined as follows:

\begin{definition}\label{D:set of sliding circuits}
We say that $y\in G$ belongs to a {\bf sliding circuit} if
$\mathfrak s^m(y)=y$ for some $m\geq 1$. Given $x\in G$, we define the
{\bf set of sliding circuits of $x$}, denoted by $\SC(x)$, as the set
of all conjugates of $x$ which belong to a sliding circuit.
\end{definition}

Since the partial order $\preccurlyeq$ is invariant under $\tau$, one has $\p(\tau(y))=\tau(\p(y))$, whence $\tau$ and $\s$ commute.  In particular, one has $y\in\SC(x)$ if and only if $y^{\Delta^k}=\tau^k(y)\in\SC(x)$ for all $k\in\ZZ$.

Our algorithm will not only compute the set $\SC(x)$, but also conjugating elements connecting the
elements of $\SC(x)$.  Basically, it constructs a connected directed graph, whose vertices
correspond to the elements of $\SC(x)$ and whose arrows correspond to conjugating elements sending
one given element in $\SC(x)$ to another.

\begin{definition}\label{D:sliding circuits graph}
Given $x\in G$, the {\bf sliding circuits graph} $\SCG(x)$ of $x$ is the directed
graph whose set of vertices is $\SC(x)$ and whose arrows correspond to conjugating elements as
follows: There is an arrow which starts at $u\in\SC(x)$, ends at $v\in\SC(x)$ and is labelled by
$s\in P\setminus\{1\}$ if and only if:
\begin{enumerate}\vspace{-\topsep}
\item $u^s=v$.

\item $s$ is an {\bf indecomposable conjugator}, that is, $s\ne 1$ and there is no element $t$,

  such that $1\prec t \prec s$ and $u^t\in \SC(x)$.
\end{enumerate}
\end{definition}

We remark that the label of each arrow is a simple element
(see Corollary~\ref{C:SCG(x) finite and connected}).

Finally, we need to define two operations that will be applied to the
conjugating elements. They are analogous to the ones defined
in~\cite{Gebhardt}, and we use the same names.

\begin{definition}\label{D:transport}
Given $x,\alpha\in G$, we define the {\bf transport} of $\alpha$ at
$x$ under cyclic sliding as
$$
 \alpha^{(1)} =\p(x)^{-1}\: \alpha \: \p(x^{\alpha}).
$$
That
is, $\alpha^{(1)}$ is the conjugating element that makes the following
diagram commutative, in the sense that the conjugating element along
any closed path is trivial:
$$
\begin{CD}
  x  @>\p(x)>>  \s(x) \\
  @V\alpha VV     @VV\alpha^{(1)}V \\
  x^{\alpha}  @>>\p(x^{\alpha})>  \s(x^{\alpha})
\end{CD}
$$
Note that the horizontal rows in this diagram correspond to
applications of cyclic sliding.

For an integer $i>1$ we define recursively
$\alpha^{(i)} = (\alpha^{(i-1)})^{(1)}$.  Note that
$(\alpha^{(i-1)})^{(1)}$ indicates the transport of $\alpha^{(i-1)}$
at $\s^{i-1}(x)$.  We also define $\alpha^{(0)} = \alpha$.
\end{definition}

The above operation is a way to transport a conjugating element
along a sliding path.  However, occasionally we will need to go backwards,
in some sense, although the obtained element will not necessarily be a
pre-image under transport.  In Section~\ref{SS:algorithm2} we will define the {\em pullback}
$s_{(1)}$ of a positive element $s$ at an element $y=\s(z)\in\SC(x)$ via the properties of its
transport at $z$ and define recursively $s_{(i)} = (s_{(i-1)})_{(1)}$ for any integer $i>1$ and
$s_{(0)}=s$ (Definition~\ref{D:pullback2}).  The details are somewhat technical and require some
prior work, so we postpone them at this stage.  At the moment, we just need to know how to compute
pullbacks in a certain special case; this is the content of the following
proposition which will be shown in Section~\ref{SS:algorithm2}:

\smallskip{\bf Proposition~\ref{P:pullback_computation}.}
{\em
Let $x\in G$, $z\in\SC(x)$, $y=\s(z)$ and let $s\in G$ be positive such that
$y^s$ is super summit.  Then the pullback of $s$ at $y$, as given in
Definition~\ref{D:pullback2}, is
$$
  s_{(1)} = \left( \p(z) \ s \ \rp(y^s)^{-1} \right) \;\vee\; 1 \;.
$$
}
\vspace{-2ex}

Hence, $s_{(1)}=\beta\vee 1$, where $\beta\in G$ is the element that makes the following
diagram commutative, in the sense that the conjugating element along
any closed path is trivial:
$$
\begin{CD}
  z @>\p(z)>>  y \\
  @V\beta VV     @VVsV \\
  \rs(y^s)  @>>\rp(y^s)>  y^s
\end{CD}
$$

\subsection{The algorithm}\label{SS:algorithm}

In this subsection we will describe in detail our algorithm to solve the CDP and the CSP in a
Garside group $G$. The only requirement needed to implement it, which we assume to be fulfilled for
the given Garside group $G$, is to know the structure of the lattices of simple elements, with
respect to both $\preccurlyeq$ and $\succcurlyeq$.
More precisely, one should have the following:
\begin{enumerate}\vspace{-\topsep}

\item A list containing the atoms,
  $\mathcal A=\{a_1,\ldots,a_\lambda\}$.

\item A function that, given $a\in \mathcal A$ and $s\in [1,\Delta]$, determines whether
$a\preccurlyeq s$ and, in that case, computes the simple element $a^{-1}s$.

\item A function that, given $a\in \mathcal A$ and $s\in [1,\Delta]$, determines whether
$s \succcurlyeq a$ and, in that case, computes the simple element  $s\:a^{-1}$.

\end{enumerate}

In Section~\ref{S:ComputingInGarsideGroups} we will see how, provided the above requirements are
fulfilled, one can compute right and left complements, gcds and lcms, normal forms, preferred
prefixes and suffixes, cyclic slidings, transports and pullbacks.

The whole algorithm is divided into three parts, called Algorithms 1, 2 and 3.  Algorithm 1
computes one element $\widetilde{x}$ in the set $\SC(x)$, starting from an arbitrary element
$x\in G$.  The algorithm also computes a conjugating element from $x$ to $\widetilde x$.
Algorithm 2 computes the arrows in the graph $\SCG(x)$ which start at a given vertex; this is
necessary for computing the entire set $\SC(x)$.
Moreover, knowing all arrows of the graph will allow us to compute a conjugating element for every
pair of elements in $\SC(x)$. Finally, Algorithm 3 solves the CDP and the CSP in $G$ using
Algorithms~1 and 2.

We remark that Algorithm 1 is a refinement of the algorithm in~\cite{EM} to compute an element in
the so-called super summit set of $x$. Here we replace two kinds of conjugation, called
{\it cycling} and {\it decycling}, by a single kind of conjugation: cyclic sliding.
This is one of the reasons that make our algorithm simpler.
Algorithm 2 is a modification of the analogous one given in~\cite{Gebhardt}, applied to cyclic 
sliding instead of cycling.
Algorithm 3 is not new, since it is implicitly or explicitly described in~\cite{EM,FG,Gebhardt} in
the context of other invariant subsets of the conjugacy class, namely super summit sets, super
summit sets with minimal simple elements, respectively ultra summit sets.  The set $\SC(x)$ is a
subset of all of these sets~\cite{GG1}.

We recommend that the reader not try to understand the algorithms at a first reading.  They
will be clarified in the following sections, where each particular step of the algorithms
will be explained in a more humane way.  See \S\ref{S:ComplexityNew} for remarks concerning efficient implementation of the algorithms.

\bigskip

\begin{minipage}{\textwidth}
\fbox{{\sc Algorithm 1:} {\bf Computing one element in $\SC(x)$}}
\begin{longtable}{lp{12.5cm}}
{\bf Input:} & $x\in G$. \\
{\bf Output:} & $\widetilde x \in \SC(x)$ and $c\in G$ such that
  $x^c=\widetilde x$.
\end{longtable}
\begin{enumerate}

\item Set $\widetilde x=x$, $c=1$ and $\mathcal T =\emptyset$.

\item While $\widetilde x\notin \mathcal T$, set
  $\mathcal T=\mathcal T \cup \{\widetilde x\}$,
  \ $c=c \cdot \p(\widetilde x)$ \ and
  \ $\widetilde x=\s(\widetilde x)$.

\item Set $y=\s(\widetilde x)$ and $d=\p(\widetilde x)$.

\item While $y\neq \widetilde x$, set $d=d \cdot \p(y)$ and $y=\s(y)$.

\item Return $\widetilde x$ and $c=c\: d^{-1}$.

\end{enumerate}
\end{minipage}
\bigskip

\begin{minipage}{\textwidth}
\fbox{{\sc Algorithm 2:} {\bf Computing the arrows in $\SCG(x)$ starting at a given vertex}}

\begin{longtable}{lp{12.5cm}}
{\bf Input:} & $v\in \SC(x)$. \\
{\bf Output:} & The set $\mathcal A_v$ of arrows in the graph $\SCG(x)$
starting at $v$.
\end{longtable}

\begin{enumerate}
\item Compute the minimal integer $N>0$ such that $\s^N(v)=v$.

\item List the atoms of $G$, say $a_1,\ldots, a_\lambda$. Set
  $\mathcal A_v=\emptyset$ and $Atoms=\emptyset$.

\item  For $t=1,\ldots, \lambda$ do:
\begin{enumerate}
\item Set $s=a_t$.
\item While $\ell(v^s)>\ell(v)$, set $s= s\cdot \left( 1 \ \vee \
  (v^s)^{-1} \Delta^{\inf(v)} \ \vee \  v^s \Delta^{-\sup(v)} \right)$.
\item If $a_t\preccurlyeq \p(v)$, then compute the iterated $N$-pullbacks
  $s,s_{(N)},s_{(2N)},\ldots$ until the first repetition, say
  $s_{(rN)}$, and set $s=s_{(rN)}$.
\item Compute the iterated $N$-transports $s,s^{(N)}, s^{(2N)},\ldots$
  until the first repetition, say $s^{(jN)}$. Let $i<j$ be such that
  $s^{(iN)}=s^{(jN)}$.
\item If $a_t\preccurlyeq s^{(mN)}$ for some $m$ with $i\leq m<j$, then do:
   \begin{enumerate}
   \item If $a_k\not\preccurlyeq s^{(mN)}$ for all $k=1,\dots,\lambda$ such that either
         $a_k\in Atoms$ or $k>t$, then set \\
         \mbox{\quad} $\mathcal A_v = \mathcal A_v \cup \{s^{(mN)}\}$
         \quad and \quad $Atoms = Atoms \cup \{a_t\}$.
   \end{enumerate}
\end{enumerate}
\item Return $\mathcal A_v$.
\end{enumerate}
\end{minipage}
\bigskip

\begin{minipage}{\textwidth}
\fbox{{\sc Algorithm 3:} {\bf Solving the conjugacy problems in $G$}}

\begin{longtable}{lp{12.5cm}}
{\bf Input:} & $x,y\in G$. \\
{\bf Output:} & - Whether $x$ and $y$ are conjugate.\newline
                - If $x$ and $y$ are conjugate, an element $c$ such
                      that $x^c=y$.
\end{longtable}

\begin{enumerate}

\item Use Algorithm 1 to compute $\widetilde x\in \SC(x)$ and
  $\widetilde y\in \SC(y)$, together with conjugating elements $c_1$
  and $c_2$ such that $x^{c_1}=\widetilde x$ and $y^{c_2}=\widetilde y$.

\item Set $\mathcal V=\{\widetilde x\}$, $\mathcal V'=\{\widetilde x\}$ and
  $c_{\widetilde x}=1$.

 \item While $\mathcal V' \neq \emptyset$, do:

 \begin{enumerate}

   \item Take $v\in \mathcal V'$.

   \item Use Algorithm 2 to compute $\mathcal A_v$.

   \item For every $s\in \mathcal A_v$, do:

   \begin{enumerate}

     \item If $v^s=\widetilde y$, then set
       $c_{\widetilde y}=c_v\cdot s$. Return `$x$ and $y$ are
       conjugate by
       $c_1 \cdot c_{\widetilde y} \cdot c_2^{-1}\;$'. STOP.

     \item If $v^s\notin \mathcal V$, then set $c_{v^s}=c_v\cdot s$,
       $\mathcal V=\mathcal V\cup \{v^s\}$, and $\mathcal V'=\mathcal V'\cup \{v^s\}$.

   \end{enumerate}

   \item Remove $v$ from $\mathcal V'$.

 \end{enumerate}

 \item Return `$x$ and $y$ are not conjugate'.

\end{enumerate}
\end{minipage}
\bigskip

\section{Cyclic sliding and the set of sliding circuits}
\label{S:cyclic sliding}

This section summarises some properties of the cyclic sliding operation, the transport map, and the
set of sliding circuits, which we require for proving the correctness of the algorithm from
Section~\ref{SS:algorithm} and for analysing its complexity.
Most of these results were obtained in~\cite{GG1} and we refer to there for further details.

\subsection*{Properties of cyclic sliding}
Cyclic sliding does not increase the canonical length.  As $G$ is of finite type, this implies that
iterated cyclic sliding starting from any $x\in G$ eventually reaches a period, that is, produces
an element of $\SC(x)$.  Moreover, iterated cyclic sliding achieves the minimal canonical length in
the conjugacy class, that is, $\SC(x)\subseteq\SSS(x)$.  More precisely, one has the following.

\begin{lemma}[{\cite[Lemma 3.4]{GG1}}]
\label{L:inf_sup_len_under_sliding}
For every $x\in G$, one has the inequalities $\inf(\s(x)) \;\geq\; \inf(x)$,
$\sup(\s(x)) \;\leq\; \sup(x)$, and $\ell(\s(x)) \;\leq\; \ell(x)$.
\end{lemma}

\begin{corollary}[{\cite[Corollary 3.5]{GG1}}]
\label{C:sliding_reaches_period}
For every $x\in G$, iterated application of cyclic sliding eventually
reaches a period, that is, there are integers $0\le i < j$ such that $\s^i(x) = \s^j(x)$.
In particular, one has $\s^k(x)\in\SC(x)$ and $\s^{j-i}(\s^k(x))=\s^{k}(x)$ for all $k\ge i$.
\end{corollary}

\begin{proposition}[{\cite[Corollary 3.9]{GG1}}]
\label{P:sliding_reaches_SSS}
For every $x\in G$, if $\ell(x)$ is not minimal in the conjugacy class of $x$, then
$\ell(x)>\ell(\s^m(x))$ for some positive integer $m < ||\Delta||$.
In particular, one has $\SC(x)\subseteq\SSS(x)$.
\end{proposition}

\subsection*{Properties of the transport map}
Under certain (mild) assumptions, the transport map respects many aspects of the Garside structure
of $G$.  In particular, transport at super summit elements preserves positive elements and powers
of $\Delta$, and it respects the partial order $\preccurlyeq$ as well as gcds with respect to
$\preccurlyeq$.  One has:

\begin{proposition}
\label{P:transport}
Let $x\in G$ and let $\alpha,\beta\in G$ such that
$x,x^\alpha,x^\beta,x^{\alpha\wedge\beta}\in\SSS(x)$ and consider transports at $x$.  Then the
following hold.
\begin{enumerate}\vspace{-\topsep}
 \item If $\alpha$ is positive then $\alpha^{(1)}$ is positive.
 \item If $\alpha$ is positive then $\p(x) \preccurlyeq  \alpha \:\p(x^{\alpha})$.
 \item If $\alpha=\Delta^k$ for $k\in\ZZ$ then $\alpha^{(1)}=\Delta^k$.
 \item If $\alpha\preccurlyeq \beta$ then $\alpha^{(1)}\preccurlyeq \beta^{(1)}$.
 \item If $\alpha$ is simple then $\alpha^{(1)}$ is simple.
 \item $(\alpha\wedge\beta)^{(1)}=\alpha^{(1)}\wedge\beta^{(1)}$.
\end{enumerate}
\end{proposition}

\begin{proof}
Claim 1 follows from \cite[Lemma 3.15]{GG1} and is equivalent to Claim 2, as
$\alpha^{(1)}=\p(x)^{-1}\alpha\p(x^\alpha)$.
Claims 3, 4, 5 and 6 are special cases of \cite[Lemma 3.16]{GG1}, \cite[Corollary 3.18]{GG1},
\cite[Corollary 3.19]{GG1} and \cite[Proposition 3.20]{GG1}, respectively.
\end{proof}

Applying iterated cyclic sliding to a conjugate $y^s$ of $y\in\SC(x)$ will eventually produce
another element of $\SC(x)$ by Corollary~\ref{C:sliding_reaches_period}.  The following Lemma makes
this more precise:  iterated transport of $s$ along the sliding circuit of $y$ eventually becomes
periodic and this happens exactly when $\SC(x)$ has been reached.

\begin{lemma}[{\cite[Lemma 3.32]{GG1}}]
\label{L:fixed_points_under_transport} Let $x\in G$, $y\in \SC(x)$ and $s\in G$ such
that $y^s\in \SSS(x)$. Let $N$ be a positive integer such that $\s^{N}(y)=y$ and for
integers $i\ge 0$ consider the transports $s^{(iN)}$ at $y$.  Then the following hold.
\begin{enumerate}\vspace{-\topsep}
\item  There are integers $i_2\ > i_1 \ge 0$ such that
       $s^{(i_1N)}=s^{(i_2N)}$.
\item  $y^s\in \SC(x)$ if and only if there is a positive integer $k$
       such that $s^{(kN)}=s$.
\end{enumerate}
\end{lemma}

\subsection*{Convexity properties and connectedness of the sliding circuits graph}
It is well known that for any $x\in G$, the set of elements conjugating $x$ to an element in
$\SSS(x)$ is closed under $\wedge$.  This has become known as {\it convexity} and in particular
implies the existence of a minimal positive element conjugating $x$ to an element in $\SSS(x)$.

\begin{proposition}[{\cite[Proposition 4.12]{FG}} or {\cite[Proposition 3.29]{GG1}}]
\label{P:sss_gcd}
Let $x\in G$. If $x^\alpha, x^\beta \in \SSS(x)$ for elements $\alpha,\beta \in G$, then
$x^{\alpha\wedge \beta} \in \SSS(x)$.
\end{proposition}

\begin{corollary}[{\cite[Theorem 2.4]{LeeLee_AbelianSubgroups}} or {\cite[Corollary 3.30]{GG1}}]
\label{C:sss_lcm}
Let $x\in G$. If $x^\alpha, x^\beta \in \SSS(x)$ for elements $\alpha,\beta \in G$, then
$x^{\alpha\vee \beta} \in \SSS(x)$.
\end{corollary}

\begin{corollary}[{\cite[Corollary 3.31]{GG1}}]
\label{C:sss_minimal} Let $x\in G$.  There is a unique positive element $\rho(x)$ (possibly
trivial) satisfying the following.
\begin{enumerate}
\vspace{-\topsep}
\item $x^{\rho(x)}\in \SSS(x)$.
\item $\rho(x)\preccurlyeq\alpha$ for every positive $\alpha\in G$ satisfying $x^\alpha\in\SSS(x)$.
\end{enumerate}
\end{corollary}

The analogous properties for $\SC(x)$ were shown in~\cite{GG1}.  They in particular imply that
$\SCG(x)$ is a finite and connected directed graph.

\begin{proposition}[{\cite[Proposition 3.33]{GG1}}]
\label{P:sc_gcd}
Let $x\in G$. If $x^\alpha, x^\beta \in \SC(x)$ for elements $\alpha,\beta \in G$, then
$x^{\alpha\wedge \beta} \in \SC(x)$.
\end{proposition}

\begin{corollary}[{\cite[Corollary 3.34]{GG1}}]
\label{C:sc_minimal}
Let $x\in G$.  There is a unique positive element $c(x)$ (possibly trivial) satisfying the following.
\begin{enumerate}
\vspace{-\topsep}
\item $x^{c(x)}\in \SC(x)$.
\item $c(x)\preccurlyeq\alpha$ for every positive $\alpha\in G$ satisfying $x^\alpha\in \SC(x)$.
\end{enumerate}
\end{corollary}

\begin{corollary}[{\cite[Corollary 3.35]{GG1}}]
\label{C:SCG(x) finite and connected}
For every $x\in G$, the graph $\SCG(x)$ is finite and connected. Moreover, the arrows of $\SCG(x)$
correspond to simple elements, and the number of arrows starting at a given vertex is bounded above
by the number of atoms of $G$.
\end{corollary}

\subsection*{Cyclic right sliding and right transport}
Recall that in a Garside group $G$ with Garside structure $(G,P,\Delta)$, apart from the prefix
order~$\preccurlyeq$, one also has the suffix order $\succcurlyeq$, defined by $a\succcurlyeq b$ if
and only if $ab^{-1}\in P$. With respect to the latter, one can consider the the notions of
preferred suffix, cyclic right sliding and set of right sliding circuits (denoted $\SC^{\Lsh}(x)$),
which are analogous to those of preferred prefix, cyclic sliding and set of sliding circuits, but
refer to the partial order $\succcurlyeq$ instead of $\preccurlyeq$ (cf.~Definitions~\ref{D:preferred prefix}, \ref{D:cyclic sliding}
and~\ref{D:set of sliding circuits}).

Consequently, one can also define a transport map for cyclic right sliding, as follows.  We remark
that, when one considers these notions with respect to $\succcurlyeq$, and tries to relate them to
the analogous notions with respect to $\preccurlyeq$, one must consider conjugating elements {\it
on the left}, meaning that a (left) conjugating element $\alpha$ relates $x$ to
$x^{\alpha^{-1}} = \alpha x \alpha^{-1}$.

\begin{definition}
\label{D:right_transport}
Given $x,\alpha\in G$, we define the {\bf right transport} of $\alpha$ at $x$ under cyclic right
sliding as $\alpha^{\rt{(1)}} = \rp(x^{\alpha^{-1}})\: \alpha \: \rp(x)^{-1}$. That
is, $\alpha^{\rt{(1)}}$ is the conjugating element that makes the following
diagram commutative, in the sense that the conjugating element along
any closed path is trivial:
\[
\xymatrix@C=20mm@R=12mm{
x
  & \rs(x) \ar[l]_(0.55){\rp(x)} \\
x^{\alpha^{-1}} \ar[u]^{\alpha}
  & \rs(x^{\alpha^{-1}}) \ar[l]^(0.55){\rp(x^{\alpha^{-1}})}
    \ar[u]_(0.55){\alpha^{\rt{(1)}}}
}
\]
\end{definition}

All results for cyclic (left) sliding and (left) transport hold in analogous form for cyclic right 
sliding and right transport; the proofs can be translated in a straight-forward way.
Alternatively, one can consider a different Garside structure.  As shown in~\cite{GG1},
$(G,P^{-1},\Delta^{-1})$ also is a Garside structure for $G$, called the {\bf reverse Garside
structure}, and cyclic right sliding and right transport with respect to $(G,P,\Delta)$ are just 
cyclic (left) sliding and (left) transport with respect to $(G,P^{-1},\Delta^{-1})$.  We refer
to~\cite[\S3.3.2]{GG1} for details.  In particular, we have the following right versions of
Lemma~\ref{L:inf_sup_len_under_sliding} and Proposition~\ref{P:transport} (1).

\begin{lemma}
\label{L:inf_sup_len_under_right_sliding}
For  $x\in G$, one has $\inf(\rs(x)) \geq \inf(x)$, $\sup(\rs(x)) \leq \sup(x)$, and
$\ell(\rs(x)) \leq \ell(x)$.
\end{lemma}

\begin{proposition}
\label{P:right_transport}
Let $x\in G$ and let $\alpha\in G$ be positive such that $x,x^{\alpha^{-1}}\in\SSS(x)$. Then, the right
transport $\alpha^{\rt{(1)}}$ of $\alpha$ at $x$ is positive.
\end{proposition}

A relation between cyclic (left) sliding and cyclic right sliding is given by the following result.

\begin{proposition}[{\cite[Proposition 3.26]{GG1}}]
\label{P:preferred_prefix_suffix}
Let $x\in G$.  Then for any $z\in \SSS(x)$ one has
$\rp(\s(z))\succcurlyeq\p(z)$ and $\rp(z)\preccurlyeq\p(\rs(z))$.
\end{proposition}

\section{Description of the algorithm}\label{S:description of the algorithm}

In this section we will explain the algorithms from \S\ref{SS:algorithm} and prove their
correctness. The main idea of these algorithms, as for the previous solutions to the conjugacy
problem given in~\cite{EM,FG,Gebhardt}, is the computation of a finite subset of the conjugacy
class, which is an invariant of the conjugacy class, together with conjugating elements connecting
each pair of elements of this subset. In our case, the finite set is $\SC(x)$, the vertex set of the connected graph $\SCG(x)$, and the conjugating elements will be paths in $\SCG(x)$.

\subsection{Algorithm 3}

We start by explaining Algorithm~3 from \S\ref{SS:algorithm}. We remark that analogues of
this algorithm, which use other sets instead of $\SC(x)$, are already given
in~\cite{EM,FG,Gebhardt}. We explain the version given in this paper which uses the
invariant $\SC(x)$.

It is clear from the definition that $\SC(x)$ is an invariant subset of the conjugacy class of $x$.
Moreover, we will see that Algorithm~1 computes, given $x\in G$, an element
$\widetilde x\in \SC(x)$, that is, $\SC(x)$ is non-empty.
Hence, two elements $x$ and $y$ are conjugate if and only if $\SC(x)=\SC(y)$ or, equivalently,
$\SC(x)\cap\SC(y)\ne\emptyset$.  Thus, knowing how to compute $\SC(x)$, starting from a given
element $x$, is sufficient to solve the conjugacy decision problem.

If we also want to solve the conjugacy search problem, that is, we want to find a conjugating
element form $x$ to $y$ in case they are conjugate, then we can do the following. Since
$\SC(x)=\SC(y)$, we just need to find an element $z\in \SC(x)$, a conjugating element $c$ form $x$
to $z$, and a conjugating element~$c_2$ from $y$ to $z$. Then $c\: c_2^{-1}$ conjugates $x$ to $y$.
In order to obtain these conjugating elements, we proceed as follows.

Suppose that $x,y\in G$ are conjugate. As we shall see, Algorithm~1 computes,  given $x\in G$, an
element $\widetilde x\in \SC(x)$ and a conjugating element $c_1$ from $x$ to $\widetilde x$.
Applying the same algorithm to~$y$, we obtain an element $\widetilde y\in \SC(y)=\SC(x)$ and a
conjugating element $c_2$ from $y$ to $\widetilde y$.  Hence, in order to obtain a conjugating
element from $x$ to $y$, we just need to find a conjugating element from $\widetilde x$ to
$\widetilde y$. In other words, we need to know how to relate, through a conjugation, any pair of
elements of $\SC(x)$. This is achieved thanks to the connected graph $\SCG(x)$, since the vertices
of this graph correspond to the elements in $\SC(x)$, and a path between two vertices corresponds
to a conjugating element from one vertex to the other.

Algorithm~3 computes a conjugating element from $\widetilde x$ to any other element in $\SC(x)$, by
computing a maximal tree of the graph $\SCG(x)$. More precisely, the algorithm starts in step~2 by
considering $\mathcal V=\mathcal V'=\{\widetilde x\}$ and $c_{\widetilde x}=1$. The set
$\mathcal V$ contains the elements which we know belong to $\SC(x)$, so at the beginning it only
contains $\widetilde x$. The set $\mathcal V'$ contains the elements of $\mathcal V$ that have not
yet been used in step~3 of the algorithm, so at the beginning $\mathcal V'=\mathcal V$.
Finally, whenever a new element $v$ is added to $\mathcal V$ (and also to $\mathcal V'$), we
compute an element $c_v$, which is a conjugating element from $\widetilde x$ to~$v$. Of course, in
step~2 of the algorithm, the conjugating element from $\widetilde x$ to $\widetilde x\in\mathcal V$
is $c_{\widetilde x}=1$.

Now step~3 does the following: For a known element of $\SC(x)$ which has not been processed before,
that is, for some $v\in \mathcal V'$, it calls Algorithm~2 to compute the arrows of $\SCG(x)$
starting at $v$. For each such arrow $s$, it computes the endpoint $v^s$ of the arrow.  If $v^s$ is
not in $\mathcal V$, this means that we encountered a new element of $\SC(x)$, so we add it to both
$\mathcal V$ and $\mathcal V'$, and at the same time compute a conjugating element from
$\widetilde x$ to $v^s$:  Since we know a conjugating element $c_{v}$ from $\widetilde x$ to $v$
and a conjugating element $s$ from $v$ to $v^s$, we can store $c_{v^s}= c_v \cdot s$ as conjugating
element from $\widetilde x$ to $v^s$. Notice that the procedure checks whether $v^s=\widetilde y$,
since in this case we have already found a conjugating element $c_{\widetilde y}$ from
$\widetilde x$ to $\widetilde y$ as desired. Concatenating it from the left with the conjugating
element from $x$ to $\widetilde x$ and from the right with the conjugating element from
$\widetilde y$ to $y$, this produces a conjugating element from $x$ to $y$ which becomes the output
of the algorithm. If $\widetilde y$ is not encountered, we remove $v$ from $\mathcal V'$ at the end
of step~3 in order to record the fact that the arrows starting at $v$ have been processed.

Notice that the procedure in step~3 is repeated while $\mathcal V'\neq \emptyset$. Since
$\mathcal V\subseteq \SC(x)$, where $\SC(x)$ is a finite set, since every element of $\mathcal V$
is added to $\mathcal V'$ exactly once, and since the procedure removes one element from
$\mathcal V'$ each time it is executed, this means that at some point we will have
$\mathcal V'=\emptyset$ and the procedure will stop. At this point, the arrows starting at every
element of $\mathcal V$ have been processed (exactly once).  Moreover, one has $\mathcal V=\SC(x)$,
since otherwise there would exist some element $v\in\mathcal V$ and some element
$w\in\SC(x)\backslash\mathcal V$ such that there is an arrow in $\SCG(x)$ from $v$ to~$w$.  (This
follows, since the graph $\SCG(x)$ is connected by \ref{C:SCG(x) finite and connected}.) But since
$v \in \mathcal V$ and $\mathcal V'=\emptyset$, step~3 has been applied to $v$, which means that
$w$ has been added to the set $\mathcal V$, a contradiction. Therefore, when the procedure stops,
one has $\mathcal V=\SC(x)$. If $\widetilde y$ was not found in $\mathcal V$, this means that
$\widetilde y\notin \SC(x)$, whence $x$ and $y$ are not conjugate.

Therefore, Algorithm~3 solves the conjugacy decision problem and the conjugacy search problem in
Garside groups of finite type, provided that Algorithms 1 and 2 are correct.

\subsection{Algorithm 1}\label{SS:algorithm1}

Given $x\in G$, Algorithm~1 finds one element $\widetilde x\in \SC(x)$ and a conjugating element
$c$ such that $x^c=\widetilde x$. This is achieved by iterated applications of cyclic sliding to
$x$.
By Corollary~\ref{C:sliding_reaches_period}, there must exist two positive integers
$0\leq i<j$ such that $\s^i(x)=\s^j(x)$, that is, $\s^i(x)\in \SC(x)$.
%
Algorithm~1 computes this element $\s^i(x)$, where $i$ is minimal. This is done by storing
all the elements $\{\s^m(x)\ |\ m\geq 0\}$, the {\em trajectory} of $x$ under cyclic sliding,
in a set called $\mathcal T$. Initially, one has
$\mathcal T=\emptyset$ and $\widetilde x=x$.  At the beginning of the $k$-th iteration of the loop
in step~2, one has $\mathcal T=\{\s^0(x), \s^1(x), \ldots, \s^{k-2}(x)\}$ and
$\widetilde x=\s^{k-1}(x)$. If $\widetilde x \notin \mathcal T$, then $\widetilde x$ is added to
$\mathcal T$ and cyclic sliding is applied to $\widetilde x$ before the next iteration of the loop.
Otherwise, a repetition (the first one) has been found and the loop terminates.

Moreover, $c$ is at every time a conjugating element from $x$ to $\widetilde x$: At the beginning
of the first iteration of the loop in step~2, $c=1$ is a conjugating element from $x$ to
$\widetilde x =x$. In each iteration of the loop, the element $c$, which is a conjugating element
from $x$ to $\widetilde x$, is multiplied on the right by $\p(\widetilde x)$, yielding a
conjugating element from $x$ to $\s(\widetilde x)$, and $\widetilde x$ is replaced by
$\s(\widetilde x)$.

Therefore, when the loop of step~2 stops, $\widetilde x = \s^i(x)\in\SC(x)$ (with $i$
minimal) and $c$ is a conjugating element from $x$ to $\widetilde x$, as desired. But notice that
the conjugating element $c$ is unnecessary long, as it contains, as a suffix, the product of all
conjugating elements along the sliding circuit containing~$\widetilde x$.  Steps~3 and 4
remove this suffix from $c$.

Step~3 initialises $y=\s(x)$ and $d=\p(x)$. The loop in step~4 checks whether $y=\widetilde x$,
otherwise applies cyclic siding to $y$ and multiplies $d$ by the corresponding conjugating element,
$\p(y)$, in such a way that when the loop terminates, the element $d$ equals the product of all
conjugating elements along the sliding circuit containing $\widetilde x$. The algorithm then
returns $\widetilde x\in\SC(x)$ and $cd^{-1}$ as the conjugating element from $x$
to $\widetilde x$.

\subsection{Algorithm 2}
\label{SS:algorithm2}

Algorithm~2 is the most involved among all the procedures in this paper. It takes an element
$v\in\SC(x)$, that is, a vertex of the graph $\SCG(x)$, and computes the arrows of $\SCG(x)$
starting at $v$. In other words, Algorithm~2 computes the indecomposable conjugators from $v$ to
other elements of $\SC(x)$. To show the correctness of each step of the algorithm, we first need to 
prove some theoretical results.

Proposition~\ref{P:sc_gcd} and Corollary~\ref{C:sc_minimal} have a consequence which is crucial for
computing the sliding circuits graph of an element $x\in G$:  Given $y\in\SC(x)$ and $s\in G$, there
is a unique $\preccurlyeq$-minimal element $c_s=c_s(y)$ satisfying $s\preccurlyeq c_s$ and
$y^{c_s}\in\SC(x)$; specifically, $c_s(y)=s\cdot c(y^s)$.
Moreover, as $s\preccurlyeq\Delta^{\sup(s)}$ and $x^{\Delta^{\sup(s)}}\in\SC(x)$, we have
$c_s(y)\preccurlyeq\Delta^{\sup(s)}$.
In particular, the number of arrows starting at a given vertex $y$ of $\SCG(x)$ is bounded by the
number of atoms of $G$ and the label of each arrow is a simple element.  In order to find the arrows
starting at $y$ it is hence sufficient to consider the set of simple elements
$\{ c_a(y)\ |\ a\textrm{ is an atom of $G$} \}$.  Let us then see how to compute $c_s(y)$ given
$y\in \SC(x)$ and $s\in G$.

By Lemma~\ref{L:fixed_points_under_transport}, the element $c_s$ we are looking for is a fixed
point under some power of transport along the sliding circuit containing $y$, and we know by
Proposition~\ref{P:transport}~(4) that transport of conjugating elements between super summit
elements respects the partial order~$\preccurlyeq$.
The basic idea is to apply iterated transport to a suitable element $p_s$, which is derived from
$s$ and satisfies $s\preccurlyeq p_s\preccurlyeq c_s$, until that fixed point is reached.  All we
need to do is to ensure that $y^{p_s}$ is super summit (so that $\preccurlyeq$ is respected) and
that $s\preccurlyeq p_s^{(kN)}$ for a sufficiently large multiple $kN$ of the length $N$ of the
sliding circuit containing $y$ (so that we can be sure that we obtain the ``right'' fixed point,
that is, one which has $s$ as a prefix).

The first step in the computation of $p_s$ is to find an element $\rho_s$ which satisfies
$s\preccurlyeq\rho_s$ and $y^{\rho_s}\in \SSS(x)$, and which is $\preccurlyeq$-minimal among all
elements doing so; this is due to~\cite{FG}.  Note that $\rho_s\preccurlyeq c_s$ since
$\SC(x)\subseteq \SSS(x)$. By Corollary~\ref{C:sss_minimal}, we have $\rho_s = s\cdot\rho(y^s)$, so
we just need to be able to compute $\rho(y^s)$.  This is achieved by the following result.

\begin{proposition}
\label{P:alg_rho_s}
For $x\in G$, the following algorithm computes $\rho(x)$ as in Corollary~\ref{C:sss_minimal}.

\begin{minipage}{\textwidth}
\textup{
\begin{enumerate}\vspace{-\topsep}
\item Set $\rho=1$.
\item While $\inf(x^{\rho})<\sinf(x)$ or $\sup(x^{\rho})>\ssup(x)$ do:
   \begin{enumerate}
   \item Set $\rho = \rho
        \cdot(1 \;\vee\; (x^{\rho})^{-1}\Delta^{\sinf(x)}
                \;\vee\; x^{\rho}\Delta^{-\ssup(x)})$.
   \end{enumerate}
\item Return $\rho(x)=\rho$.
\end{enumerate}
}
\end{minipage}

If $x=y^s$ with $y\in\SSS(x)$, then the algorithm terminates after at most
$\ell(s)\cdot\!\!\wordlen{\Delta}$ passes through the loop.
\end{proposition}
\begin{proof}
Let $\alpha$ be a positive element such that $x^{\rho\alpha}\in\SSS(x)$.
Then $(x^{\rho\alpha})^{-1}\in\SSS(x^{-1})$ and we have $\sup((x^{\rho\alpha})^{-1})=-\sinf(x)$ and
$\sup(x^{\rho\alpha})=\ssup(x)$.  Thus $x^{\rho}\preccurlyeq x^{\rho}\alpha = \alpha x^{\rho\alpha}
\preccurlyeq\alpha\Delta^{\ssup(x)}$, whence $x^{\rho}\Delta^{-\ssup(x)}\preccurlyeq\alpha$ and,
analogously, $(x^{\rho})^{-1}\Delta^{\sinf(x)}\preccurlyeq\alpha$.
As $\alpha$ is positive, the above implies $1 \;\vee\; (x^{\rho})^{-1}\Delta^{\sinf(x)}
\;\vee\; x^{\rho}\Delta^{-\ssup(x)}\preccurlyeq\alpha$. Moreover,
$1\;\vee\;(x^{\rho})^{-1}\Delta^{\sinf(x)}\;\vee\; x^{\rho}\Delta^{-\ssup(x)}=1$ if and only if
$\sup(x^{\rho})\le\ssup(x)$ and $\inf(x^{\rho})=\sup((x^{\rho})^{-1})\ge\sinf(x)$, that is, if
and only if $x^{\rho}\in\SSS(x)$.

Hence, at any stage of the above algorithm, the element $\rho$ satisfies $\rho\preccurlyeq c$ for
every positive element $c\in G$ such that $x^c\in\SSS(x)$.  In particular, $\wordlen{\rho}$ is
bounded. As $\wordlen{\rho}$ is strictly increasing at every step of the algorithm, the algorithm
terminates and outputs $\rho(x)$ as claimed. Finally, if $x=y^s$ with $y\in\SSS(x)$, then
$\rho(x)\preccurlyeq s^{-1}\Delta^{\sup(s)}\preccurlyeq\Delta^{\ell(s)}$, whence the algorithm
terminates after at most $\ell(s)\cdot\!\!\wordlen{\Delta}$ steps.
\end{proof}

\begin{corollary}
\label{C:alg2_step3b}
Steps~3\,(a) and 3\,(b) in Algorithm~2 compute the element $\rho_{a_t}$. The body of the while loop
is executed at most $\wordlen{\Delta}$ times.
\end{corollary}

\begin{proof}
Note that in steps~3\,(a) and 3\,(b) of Algorithm~2 we have $v\in\SC(x)\subseteq\SSS(x)$ and
$\ell(a_t)=1$.  Hence, by Proposition~\ref{P:alg_rho_s}, steps~3\,(a) and 3\,(b) in Algorithm~2
compute exactly $a_t\cdot\rho(v^{a_t})=\rho_{a_t}$ with at most $\wordlen{\Delta}$ passes through
the while loop.
\end{proof}

The element  $1 \;\vee\; (x^{\rho})^{-1}\Delta^{\sinf(x)}\;\vee\; x^{\rho}\Delta^{-\ssup(x)}
= \big(1 \;\vee\; (x^{-1})^{\rho}\Delta^{-\ssup(x^{-1})}\big)\;\vee\; \big(1 \;\vee\;
x^{\rho}\Delta^{-\ssup(x)}\big)$ in step~2\,(a) of the algorithm in
Proposition~\ref{P:alg_rho_s} can be computed efficiently using the following result.

\begin{proposition}
\label{P:alg_rho_s_lcm}
If $x\in G$ such that $\sup(x)=q+r$ with $0\leq r \leq \ell(x)$, then $1\vee x\Delta^{-q}$ is the
product of the leftmost $r$ factors of the right normal form of $x$.
\end{proposition}

\begin{proof}
Observe that $a\preccurlyeq b$ is equivalent to $b^{-1}\succcurlyeq a^{-1}$ for all $a,b\in G$ by
the definitions of $\preccurlyeq$ and $\succcurlyeq$.  This implies that one has
$a\vee b = (a^{-1}\rwedge b^{-1})^{-1}$ for all $a,b\in G$.
Hence, in particular, $1\vee x\Delta^{-q} = (1\rwedge \Delta^q x^{-1})^{-1}
 = \big((x\rwedge \Delta^q)x^{-1}\big)^{-1} = x(x\rwedge \Delta^q)^{-1}$.  Since
$x\rwedge \Delta^q$ contains all but the leftmost $r$ factors of
the right normal form of $x$, the claim follows.
\end{proof}

Next we consider the sequence of iterated transports along the sliding circuit which contains the
element $y$.  This sequence will eventually become periodic; we are interested in the periodic part.

\begin{definition}
\label{D:fixed_points}
Let $x\in G$, $y\in \SC(x)$ and $u\in G$ such that $y^u\in \SSS(x)$ and let $N$ be the length of
the sliding circuit containing $y$, that is, let $N$ be the smallest positive integer such that
$\s^N(y)=y$.  For integers $i\ge 0$ consider the transports $u^{(iN)}$ at $y$. By
Lemma~\ref{L:fixed_points_under_transport}, there are integers $i_2 > i_1 \ge 0$ such that
$u^{(i_1N)} = u^{(i_2N)}$. Let $i_1$ and $i_2$ be minimal subject to this condition and define
$l(u) = i_2 - i_1$ and $F(u) = \{ u^{(iN)}\ |\ i_1\le i < i_2 \}$.
\end{definition}

Note that $1 \in F(u)$ if and only if $F(u) = \{ 1 \}$.  Moreover by
Lemma~\ref{L:fixed_points_under_transport}, $y^u\in \SC(x)$ if and only if $i_1=0$ , that is, if
and only if $u\in F(u)$. Finally, for all $v\in F(u)$ and all $i\in\NN$ we have $v^{(il(u)N)}=v$,
in particular, $y^v\in \SC(x)$.  In other words, the set $F(u)$ contains those iterated transports
$u^{(iN)}$ of~$u$ along the sliding circuit of $y$ which are fixed by some iterated transport along the sliding circuit, that is, those which satisfy $y^{u^{(iN)}}\in\SC(x)$.

\begin{lemma}
\label{L:no_proper_prefix}
Let $x\in G$, $y\in \SC(x)$, $s \in G$ and denote $c_s=c_s(y)$. Let $N$ be the length of the sliding
circuit containing $y$, that is, let $N$ be the smallest positive integer such that $\s^N(y)=y$. 
If $c_s \preccurlyeq c_s^{(iN)}$ for some $i>0$ then $c_s^{(iN)} = c_s$.
\end{lemma}

\begin{proof}
Let $c_s^{(iN)} = c_s \gamma$ with a positive element $\gamma$. By induction,
$c_s\gamma\preccurlyeq c_s^{(kiN)}$ for all $k\ge 1$ from Proposition~\ref{P:transport}~(4).  By
the preceding remark, this in particular implies
$c_s \preccurlyeq c_s \gamma \preccurlyeq c_s^{(l(c_s) iN)} = c_s$, that is, $\gamma = 1$.
\end{proof}

\begin{lemma}
\label{L:fixed_points_of_transport}
Let $x\in G$, $y\in \SC(x)$, $s \in P$ and denote $c_s=c_s(y)$. Assume that $p\in P$ satisfies
$p\preccurlyeq c_s$ and $y^p \in\SSS(y)$ and that $F = F(p)\neq \{1\}$.
\begin{enumerate}
\vspace{-\topsep}
\item If there exists $v \in F$ such that $s\preccurlyeq v$ then
      $c_s = v$.
\item If $s\not\preccurlyeq v$ for all $v\in F$, then
      $c_s$ is not an indecomposable conjugator starting at $y$.
\end{enumerate}
\end{lemma}

\begin{proof}
First note that by Proposition~\ref{P:transport}~(4), we have $p^{(i)} \preccurlyeq c_s^{(i)}$ for
all $i>0$.

Assume first that $v \in F$ such that $s\preccurlyeq v$.  Then $y^v\in \SC(x)$ and the minimality
of $c_s$ implies $c_s\preccurlyeq v$.  Now $v = p^{(iN)}$ for some $i$, whence
$c_s\preccurlyeq v = p^{(iN)} \preccurlyeq c_s^{(iN)}$. Lemma~\ref{L:no_proper_prefix} then yields
$v = c_s$ and Claim 1 is shown.

Now assume that $s\not\preccurlyeq v$ for all $v\in F$ and let $i$ be a multiple of $l(c_s)$
sufficiently large so that $v = p^{(iN)} \in F$. Since $1\notin F$, we have $v\ne 1$ and
$y^v\in\SC(x)$. Moreover, $v = p^{(iN)} \preccurlyeq c_s^{(iN)} = c_s$ and $v\ne c_s$, since
$s\not\preccurlyeq v$ but $s\preccurlyeq c_s$.  Hence, $c_s$ is not an indecomposable conjugator starting at $y$ and Claim~2 is shown.
\end{proof}

Recall that we are trying to compute the arrows of $\SCG(x)$ starting at $y$. In Algorithm~2, we
start with an atom $a$ and we try to see if there is an arrow $c$ starting at $y$ such that
$a\preccurlyeq c$ or, equivalently, such that $\rho_a\preccurlyeq c$. The above result says that if
$F(\rho_a)\neq \{1\}$ then we will have no problem, since either $c_a$ can be computed by iterated
transport (where $c_a$ is the only possible candidate for being such an arrow), or we can be sure
that $c_a$ is not an arrow, since it is decomposable.  Unfortunately, it may occur that
$F(\rho_a)=\{1\}$, as we can see in the following example:

\begin{example}
\label{ex_trivial_transport}
Consider in the Artin braid group $B_5$ the elements
$y=x=\Delta\cdot\sigma_2\sigma_1\sigma_4\sigma_3\sigma_4\cdot\sigma_1$, in left normal form as
written, and $s = \sigma_3\sigma_2\sigma_1$. It is easy to check that $\s^6(y) = y$, that is,
$y\in\SC(x)$. Since $y^s=\Delta\cdot\sigma_1\sigma_3\cdot\sigma_3\sigma_2\sigma_1\sigma_2$ is in
left normal form as written, $y^s\in \SSS(x)$, that is, $\rho_s=s$.

However, $s^{(1)} = \p(y)^{-1}s\p(y^s) = 1$ and hence $F(s) = \{1\}$, that is, the requirements of
Lemma~\ref{L:fixed_points_of_transport} are not satisfied.
\end{example}

The above example shows that one could possibly have $F(\rho_a)=\{1\}$ for some atom $a$ in the
situation of Algorithm~2.  In this case, Lemma~\ref{L:fixed_points_of_transport} would not
guarantee that iterated transport is sufficient to find $c_a$ or to be sure that $c_a$ is
decomposable.  Let us now see that there is another condition which also ensures that either $c_a$
can be computed by iterated transport, or that it is decomposable; it is given by the corollary to
the following result.

\begin{lemma}
Let $x\in G$ and $v\in \SC(x)$. Let $s\neq 1$ be a positive element such that $v^s\in \SSS(x)$. If
$s^{(k)}=1$ for some $k\geq 1$, then $s\wedge \p(v)\neq 1$.
\end{lemma}

\begin{proof}
This proof parallels the one of \cite[Lemma~4.11]{Gebhardt}. Denote $w=v^s$. By hypothesis
$$
  s^{(k)}=\big(\p(v) \p(\s(v))\cdots\p(\s^{k-1}(v))\big)^{-1}
  \: s \: \big(\p(w) \p(\s(w))\cdots \p(\s^{k-1}(w))\big)=1,
$$
that is,
$$
 s \: \big(\p(w) \p(\s(w))\cdots \p(\s^{k-1}(w))\big) = \p(v) \p(\s(v))\cdots \p(\s^{k-1}(v)).
$$
We will show the result by induction on $k$. If $k=1$, one has $s\p(w)= \p(v)$, hence
$s\wedge\p(v)=s \neq 1$.  Suppose the result is true for $k-1$, and consider $s^{(1)}$. We can
assume that $s^{(1)}\neq 1$, otherwise the result would hold by applying the case $k=1$. But we
have $(s^{(1)})^{(k-1)}=1$, so by induction hypothesis $s^{(1)}\wedge \p(\s(v))\neq 1$.

Recall that the transport $t^{(1)}$ of an element $t$ at $v$ satisfies
$t^{(1)}=\p(v)^{-1}t\p(v^t)$.  For $t=\p(v)$ this yields $\p(v)^{(1)} = \p(v^{\p(v)}) = \p(\s(v))$.
As the transport preserves $\wedge$ by Proposition~\ref{P:transport}~(6), one hence has
$(s \wedge \p(v))^{(1)} = s^{(1)}\wedge \p(v)^{(1)} = s^{(1)}\wedge \p(\s(v))\neq 1$, which implies
$s\wedge \p(v)\neq 1$ by Proposition~\ref{P:transport}~(3).
\end{proof}

\begin{corollary}
\label{C:fixed_points_of_transport}
Let $x\in G$ and $v\in \SC(x)$. Let $a$ be an atom such that $a\not\preccurlyeq \p(v)$. Then either
$F(\rho_a)\neq \{1\}$ or $c_a$ is not an indecomposable conjugator starting at $v$.
\end{corollary}

\begin{proof}
Suppose that $F(\rho_a)=\{1\}$. This means that some iterated transport $(\rho_a)^{(k)}=1$ for some
$k\geq 1$. By the above lemma, it follows that $\rho_a\wedge \p(v)\neq 1$. Hence there must exist
an atom $b$ such that $b\preccurlyeq \rho_a \wedge \p(v)$. Since $b\preccurlyeq \p(v)$ and
$v^{\p(v)}\in \SC(x)$, it follows that $c_b\preccurlyeq \p(v)$. On the other hand, since
$b\preccurlyeq \rho_a \preccurlyeq c_a$, it follows that $c_b\preccurlyeq c_a$. But one cannot have
$c_b=c_a$, otherwise $a\preccurlyeq c_a=c_b\preccurlyeq \p(v)$, which is not possible by
hypothesis. Therefore, $c_b$ is a proper prefix of $c_a$, which means that $c_a$ is not an
indecomposable conjugator starting at $v$.
\end{proof}

Recall that if $F(\rho_a)\neq \{1\}$ then either $c_a$ can be found by iterated transport or $c_a$
is not indecomposable. Hence, if $a\not\preccurlyeq \p(v)$, we just need iterated transport in
order to compute or to discard $c_a$.  The case that remains to be dealt with is the case
$a\preccurlyeq \p(v)$ and $F(\rho_a)=\{1\}$.

We will now consider the more general situation that $F(\rho_s)=\{1\}$ for some element $s\in G$.
Iterated transport of $\rho_s$ reaches the ``wrong'' fixed point in this situation. The solution is
to apply iterated transport not to $\rho_s$ itself, but to a related element $p$ satisfying $\rho_s\preccurlyeq p\preccurlyeq c_s$
for which the existence of $v \in F(p)$ with $s\preccurlyeq v$ is guaranteed.  To this end we
introduce the notion of the ``pullback'' of an element $s$, defined as the $\preccurlyeq$-minimal
among the elements whose transport has $s$ as a prefix.

\begin{definition}\label{D:pullback2}
Let $x\in G$, $z\in\SC(x)$, $y=\s(z)$ and let $s\in G$ be positive. By Propositions~\ref{P:sss_gcd}
and Proposition~\ref{P:transport}~(6), there exists a unique positive $\preccurlyeq$-minimal
element $s_{(1)}\in G$ satisfying $z^{s_{(1)}} \in \SSS(x)$ and $s \preccurlyeq (s_{(1)})^{(1)}$,
where $(s_{(1)})^{(1)}$ indicates the transport of $s_{(1)}$ at $z$.  We call $s_{(1)}$ the {\bf
pullback} of $s$ at $y$.

For any integer $k>1$ we define recursively the $k$-fold pullback $s_{(k)} = (s_{(k-1)})_{(1)}$ of $s$ at $y$.  Note that $(s_{(k-1)})_{(1)}$ indicates the pullback of $s_{(k-1)}$ at the unique element $w$ in the sliding circuit of $y$ satisfying $\s^{k-1}(w)=y$.  We also define $s_{(0)} = s$.
\end{definition}

\begin{lemma}
\label{L:iterated_pullback}

Let $x\in G$, $z\in\SC(x)$, $y=\s^k(z)$ for a positive integer $k$ and let $s\in G$ be positive.
Then, the $k$-fold pullback $s_{(k)}$ of $s$ at $y$ is the minimal positive element such that
$s\preccurlyeq (s_{(k)})^{(k)}$ and $z^{s_{(k)}} \in \SSS(x)$.
\end{lemma}

\begin{proof}
The claim holds for $k=1$ by definition of the pullback.
Suppose the claim is true for $k-1$. By Proposition~\ref{P:transport}~(4), one then has
$s \preccurlyeq (s_{(k-1)})^{(k-1)} \preccurlyeq (((s_{(k-1)})_{(1)})^{(1)})^{(k-1)}
 = (s_{(k)})^{(k)}$. Moreover, if $\alpha$ is a positive element such that
$s\preccurlyeq \alpha^{(k)}$ and $z^\alpha\in\SSS(x)$, then $\alpha^{(1)}$ is a positive element
such that $s\preccurlyeq (\alpha^{(1)})^{(k-1)}$ and $\s(z)^{\alpha^{(1)}}=\s(z^\alpha)\in\SSS(x)$.
Hence, $s_{(k-1)}\preccurlyeq \alpha^{(1)}$ by induction. By definition of the pullback of
$s_{(k-1)}$, we then have $s_{(k)}=(s_{(k-1)})_{(1)}\preccurlyeq \alpha$, as we wanted to show.
\end{proof}

\begin{lemma}
\label{L:pullback_monotonic}
Let $x\in G$, $z\in\SC(x)$, $y=\s^k(z)$ for a positive integer $k$ and let  $s,t\in G$ such that
$1\preccurlyeq s\preccurlyeq t$.  Then, $s_{(k)} \preccurlyeq t_{(k)}$.
\end{lemma}

\begin{proof}
By Lemma~\ref{L:iterated_pullback}, we have $t\preccurlyeq (t_{(k)})^{(k)}$ and
$z^{t_{(k)}}\in\SSS(x)$.  But then $s\preccurlyeq t\preccurlyeq (t_{(k)})^{(k)}$ and, again using
Lemma~\ref{L:iterated_pullback}, we obtain $s_{(k)} \preccurlyeq t_{(k)}$ as we wanted to show.
\end{proof}

\begin{lemma}
\label{L:pullback_bounded}
Let $x\in G$, $z\in\SC(x)$, $y=\s(z)$ and let $s\in G$ be positive.
Then the pullback $s_{(1)}$ of $s$ at $y$ satisfies
$s_{(1)} \preccurlyeq \Delta^{\sup(s)}$.
\end{lemma}

\begin{proof}
Let $q=\sup(s)\ge 0$ and consider transport at $z$.  We have
$s\preccurlyeq\Delta^q=(\Delta^q)^{(1)}$ by Proposition~\ref{P:transport}~(3).  Moreover,
$\Delta^q$ is positive and $y^{\Delta^q}\in\SSS(x)$.  By $\preccurlyeq$-minimality of $s_{(1)}$,
we obtain $s_{(1)} \preccurlyeq \Delta^q$ as claimed.
\end{proof}

The next result shows how one can use pullbacks to compute $c_s$ in the case in which
$F(\rho_s)=\{1\}$ may occur.

\begin{proposition}
\label{P:pullback}
Let $x\in G$, $v\in \SC(x)$ and let $N$ be the length of the sliding
circuit of $v$, that is, let $N$ be the smallest positive integer such
that $\s^N(v)=v$. Let $s\in P\setminus\{1\}$ such that $v^s\in \SSS(x)$ and for
integers $k\geq 0$ consider the iterated pullbacks $s_{(kN)}$ at $v$.
Let $i\geq 0$ be such that $s_{(iN)}=s_{(jN)}$ for some $j>i$. Then
$c_s$ is the only element in $F(s_{(iN)})$ which admits $s$ as a
prefix. In particular, $F(s_{(iN)})\neq \{1\}$.
\end{proposition}

\begin{proof}
First note that by Lemma~\ref{L:pullback_bounded}, we have
$1\preccurlyeq s_{(kN)}\preccurlyeq\Delta^{\sup(s)}$ for all $k\geq 0$.  As
$G$ is of finite type, the number of such elements is finite, whence there
exist integers $i\geq 0$ and $j>i$ such that $s_{(iN)}=s_{(jN)}$.

Let $m=i(j-i)\geq i$ and denote $p=s_{(mN)}$. Notice that iterated
$N$-fold pullback becomes periodic of period $j-i$ starting from the $i$-th
term, hence  $p_{(k(j-i)N)}=p$ for all $k\geq 0$, that is, $p=s_{(k(j-i)N)}$
for all $k\geq i$.  Now recall from
Lemma~\ref{L:fixed_points_under_transport} that, since $v^{c_s}\in \SC(x)$,
we have $(c_s)^{(tN)}=c_s$ for some $t\geq 1$. Consider then $M>i$ to be a
multiple of $t$, big enough so that $p^{(M(j-i)N)}\in F(p)$. By Lemma~\ref{L:iterated_pullback},
$p=s_{(M(j-i)N)}$ is the minimal positive element such that
$s\preccurlyeq p^{(M(j-i)N)}$. This implies that $F(p)\neq\{1\}$ and that $F(p)$
contains an element which admits $s$ as a prefix. Moreover,
$s\preccurlyeq c_s = (c_s)^{(M(j-i)N)}$, where the equality in the
last step holds since $M$ is a multiple of $t$.  By the minimality of
$p$ one finally has $p\preccurlyeq c_s$.  We can then apply
Lemma~\ref{L:fixed_points_of_transport} to $p$, and conclude that
$c_s = p^{(M(j-i)N)} \in F(p)$.
Uniqueness also follows from Lemma~\ref{L:fixed_points_of_transport}.

It only remains to be shown that $F(p)=F(s_{(iN)})$, that is
$F(s_{(mN)})=F(s_{(iN)})$ for $m$ as above; indeed, we will show that
$F(s_{(kN)})=F(s_{(iN)})$ for all $k \ge i$. Since iterated $N$-fold pullback
is periodic of period $j-i$ from the $i$-th term, we can assume $i< k<j$.

We have $s_{(iN)}\preccurlyeq (s_{(kN)})^{((k-i)N)}$ and also
$s_{(kN)}\preccurlyeq (s_{(jN)})^{((j-k)N)}=(s_{(iN)})^{((j-k)N)}$ by
Lemma~\ref{L:iterated_pullback}.
Applying $(k-i)N$-fold transport to the second expression and using
Proposition~\ref{P:transport}~(4), one obtains
$(s_{(kN)})^{((k-i)N)}\preccurlyeq (s_{(iN)})^{((j-i)N)}$, whence
$s_{(iN)}\preccurlyeq (s_{(kN)})^{((k-i)N)}\preccurlyeq (s_{(iN)})^{((j-i)N)}$.

Using Proposition~\ref{P:transport}~(4) again, we can for any $K\ge 0$ apply $K$-fold transport to
this expression and we see that
$
  (s_{(iN)})^{(K)}\preccurlyeq (s_{(kN)})^{(K+(k-i)N)}
                  \preccurlyeq (s_{(iN)})^{(K+(j-i)N)}
$
for all $K\geq 0$.
That is, for any integer $K$ large enough so that $s'=(s_{(iN)})^{(K)}\in F(s_{(iN)})$,
we have $c_{s'} = s' \preccurlyeq (s')^{((j-i)N)}$ and hence $s' = (s')^{((j-i)N)}$ by
Lemma~\ref{L:no_proper_prefix} (where $s'=c_{s'}$ is chosen as the element $s$ in the
statement of the lemma).  Hence, the above inequality implies $s'=(s_{(kN)})^{(K+(k-i)N)}$.
As this is true for all sufficiently large $K$, we have $F(s_{(iN)})=F(s_{(kN)})$.  In particular,
$F(p)=F(s_{(iN)})$, whence $c_s\in F(s_{(iN)})$, as we wanted to show.
\end{proof}

The following result allows us to compute pullbacks in the situation of Algorithm~2.

\begin{proposition}
\label{P:pullback_computation}
Let $x\in G$, $z\in\SC(x)$, $y=\s(z)$ and let $s\in G$ be positive such that
$y^s$ is super summit.  Then the pullback of $s$ at $y$, as given in
Definition~\ref{D:pullback2}, is
$$
  s_{(1)} = \left( \p(z) \ s \ \rp(y^s)^{-1} \right) \;\vee\; 1 \;.
$$
\end{proposition}

\begin{proof}
Let $u=\left( \p(z)\ s \ \rp(y^s)^{-1}\right) \;\vee\; 1$. We show
that $u$ satisfies the defining properties of $s_{(1)}$.  The
following commutative diagram illustrates the situation; all
conjugating elements corresponding to arrows will be shown to be
positive.
\[
\xymatrix@C=15mm@R=12mm{
\rs(y^s)
  \ar[r]
  \ar@{.>}@<7.0ex>[rrr]^{\p(\rs(y^s))}
  \ar@{.>}@<-4.0ex>[rr]_(0.3){\rp(y^s)}
& y^{s\alpha^{-1}}
  \ar[r]^{\alpha}
  \ar@{.>}@<3.5ex>[rrr]^(0.7){\p(y^{s\alpha^{-1}})}
& y^s
  \ar[r]
  \ar@{.>}@<-4.0ex>[rrr]_(0.7){\p(y^s)}
& \s(\rs(y^s))
  \ar[r]
& \s(y^{s\alpha^{-1}})
  \ar[r]
& \s(y^s)
\\
& z
  \ar[u]_(0.35){u}
  \ar[r]_{\p(z)}
& y
  \ar[u]^(0.35){s}
  \ar[urr]_(0.35){u^{(1)}}
}
\]

\textit{Claim 1:}  $z^u \in \SSS(x)$.

 Proof:  As $y^s\in \SSS(x)$, we have $z^{\p(z)\ s \ \rp(y^s)^{-1}} = \rs(y^s)\in \SSS(x)$ by
 Lemma~\ref{L:inf_sup_len_under_right_sliding}.
 Then, Corollary~\ref{C:sss_lcm} implies $z^u \in \SSS(x)$, since
 $u=\left( \p(z)\ s \ \rp(y^s)^{-1}\right) \;\vee\; 1$.

\textit{Claim 2:}  $u$ is a positive element and $s\preccurlyeq u^{(1)}$.

 Proof: By definition, $u$ is positive and $u=\p(z) s \alpha^{-1}$,
 where $\alpha = \rp(y^s) \,\rwedge \p(z) s$.
 (Note that $\alpha^{-1} = \rp(y^s)^{-1} \,\vee\, (\p(z) s)^{-1}$.)
 By Proposition~\ref{P:preferred_prefix_suffix},
 $\rp(y^s) \preccurlyeq \p(\rs(y^s))$ and since $\rp(y^s)\alpha^{-1}$
 is positive and $\rs(y^s)\in\SSS(x)$, we obtain with
 Proposition~\ref{P:transport} (2)
 \[
 \rp(y^s) \preccurlyeq \p(\rs(y^s))
 \preccurlyeq \rp(y^s)\alpha^{-1}\p\left(
    \rs(y^s)^{\rp(y^s)\alpha^{-1}}\right)
 = \rp(y^s)\alpha^{-1}\p(y^{s\alpha^{-1}})
 \]
 which implies $1\preccurlyeq \alpha^{-1}\p(y^{s\alpha^{-1}})$.
 Hence, $s\preccurlyeq s\alpha^{-1}\p(y^{s\alpha^{-1}}) = \p(z)^{-1} u
 \p(z^u) = u^{(1)}$.

\textit{Claim 3:}  If $t$ is a positive element such that
 $z^t\in \SSS(x)$ and $s\preccurlyeq t^{(1)}$, then
 $u\preccurlyeq t$.

 Proof: Write $t^{(1)}=s\gamma$ for some positive element $\gamma$
 and apply cyclic right sliding to $y$, $y^s$ and $y^{t^{(1)}}=\s(z^t)$, as shown in the following
 commutative diagram.
 \[
 \xymatrix@C=15mm{
 z^t
   \ar[r]^{\p(z^t)}
 & \s(z^t)
 & \rs(\s(z^t))
   \ar[l]_{\rp(\s(z^t))}
 \\
 & y^s
   \ar[u]^{\gamma}
 & \rs(y^s)
   \ar[u]_{\gamma^{\rt{(1)}}}
   \ar[l]_{\rp(y^s)}
 \\
 z
   \ar[r]^{\p(z)}
   \ar[uu]^{t}
 & y
   \ar[u]^{s}
 & \rs(y)
   \ar[u]_{s^{\rt{(1)}}}
   \ar[l]_{\rp(y)}
 }
 \]
 We obtain
 $
  t = \p(z)\,s\,\gamma\,\p(z^t)^{-1}
    = \p(z)\,s\,\rp(y^s)^{-1}\,\gamma^{\rt{(1)}}\,
                \left[\rp(\s(z^t))\,\p(z^t)^{-1}\right]
 $,
 where $\gamma^{\rt{(1)}}$ is positive by Proposition~\ref{P:right_transport}
 and the factor in brackets is positive by Proposition~\ref{P:preferred_prefix_suffix}.
 Therefore, we have $\p(z)\ s\ \rp(y^s)^{-1} \preccurlyeq t$, and since
 $t$ is positive, one finally has
 $u=\left( \p(z)\ s \ \rp(y^s)^{-1}\right) \;\vee\; 1\preccurlyeq t$.
\end{proof}

\pagebreak

\begin{example}
\label{ex_pullback}
Consider the situation from Example~\ref{ex_trivial_transport}. The
trajectory of $y = \Delta \cdot \sigma_2\sigma_1\sigma_4\sigma_3\sigma_4
\cdot \sigma_1$ under cyclic sliding has length $N=6$.  Computing iterated
pullbacks of $\rho_s = s =\sigma_3 \sigma_2 \sigma_1$ at $y$ we obtain
$s_{(12)} = s_{(6)} = \sigma_3 \sigma_4$. Hence, using the notation from
Proposition~\ref{P:pullback}, we have $i=1$ and $j=2$.

Computing iterated transports of $p = s_{(iN)} = s_{(6)} = \sigma_3
\sigma_4$, we obtain
$p^{(12)} = p^{(6)} = \sigma_3 \sigma_2 \sigma_1 \sigma_4$.
Hence, we have $F(p) = \{p^{(6)}\}$ and as
$s\preccurlyeq p^{(6)}$, we obtain
$c_s = p^{(6)} = \sigma_3 \sigma_2 \sigma_1 \sigma_4$.

Note that $p\notin F(p)$, that is, computing iterated transports is
necessary even after reaching a stable loop under iterated $N$-fold pullback.
\end{example}

The results obtained in this section ensure that step~3\,(c) of Algorithm~2, if executed, will
compute an element $(\rho_{a_t})_{(iN)}$, one of whose iterated transports is
precisely $c_{a_t}$. This computation is only done whenever
$a_t\preccurlyeq \p(v)$, which is the only case, as we saw above, in which we
cannot be sure to find $c_{a_t}$ or to be able to discard $c_{a_t}$ as
decomposable using $F(\rho_{a_t})$.
Note, in particular, that computing pullbacks is not necessary if $v$ is rigid
(or, by~\cite[Theorem~1.1]{GG1} equivalently, has a rigid conjugate).
The algorithm continues in step~3\,(d) by applying
iterated transport to the corresponding element (either $\rho_{a_t}$ or
$(\rho_{a_t})_{(iN)}$) until the first repetition occurs. Then, step~3\,(e) checks
whether any of the elements in $F(\rho_{a_t})$ respectively $F((\rho_{a_t})_{(iN)})$
admits ${a_t}$ as a prefix, in which case it will precisely be
$c_{a_t}$ by Lemma~\ref{L:fixed_points_of_transport}.  If $a_t$ does not
occur as a prefix, then $c_{a_t}$ is not indecomposable by
Lemma~\ref{L:fixed_points_of_transport},
Corollary~\ref{C:fixed_points_of_transport} and
Proposition~\ref{P:pullback}.

However, even if $c_{a_t}$ occurs as an element of $F(\rho_{a_t})$ respectively
$F((\rho_{a_t})_{(iN)})$, it is not necessarily an indecomposable conjugator.
The latter property is checked in step~3\,(e)\,i: The set $Atoms$ will eventually contain the
atoms $a_k$ such that $c_{a_k}$ is an indecomposable conjugator starting at $v$ and
$k=\max\{i \ | \  a_i\preccurlyeq c_{a_k} \}$. Suppose that we
have computed $c_{a_t}$ for some atom $a_t$. If $t$ is not the biggest index
among the atoms dividing $c_{a_t}$, then we can discard $c_{a_t}$ at this step
since, if it is indecomposable, it will appear again in a further step of the
algorithm, when the mentioned atom is processed. On the other hand, if $t$ is
the maximal index among the atoms dividing $c_{a_t}$ but $c_{a_t}$ is 
decomposable, then there must exist some indecomposable $c_{a_l}\preccurlyeq
c_{a_t}$, where $l<t$ is maximal among the atoms dividing $c_{a_l}$. In particular,
$a_l$ has been processed before $a_t$, and we must have $a_l\in Atoms$.
Therefore, if $a_k\not\preccurlyeq c_{a_t}$ for all $a_k\in Atoms$
and also for all $k>t$, we can be sure that $c_{a_t}$ is
indecomposable, and we can add $a_t$ to the set $Atoms$.  This is what
is done in step~3\,(e)\,i, hence Algorithm~2 computes the arrows starting
at $v$, as claimed.

\section{Complexity of the algorithms}\label{S:complexity}

\subsection{Computing in Garside groups}\label{S:ComputingInGarsideGroups}

In this section, we will describe how one can perform all the computations required by
our algorithms in any Garside group of finite type, provided some basic operations on simple
elements can be performed. We refer the reader to~\cite{GAP} for a similar approach.

We remark that in a particular Garside group there may be specific
algorithms having better complexity than the generic ones we describe below.
This is in particular the case for braid groups (see \cite{Epstein} and \cite{BKL1}).
Hence one should not use the algorithms below if one just needs to make
computations in braid groups.

\pagebreak

\begin{assumption}
 \label{A:BasicOperations}
 Let $G$ be a Garside group of finite type.  We assume that the list of atoms
 $\mathcal A=\{a_1,\ldots,a_{\lambda}\}$ of $G$ is known and that the following operations can be
 performed effectively;
 we consider the cost of these operations to be $O(C)$.
 \begin{namedenumerate}[Op]\vspace{-\topsep}
  \item[\textup{(H)}] Given  a simple element $s$, compute a hash value for $s$.
  \item[\textup{(Op)}] Given an atom $a\in\mathcal A$ and a simple element $s$, test whether
     $a\preccurlyeq s$ (respectively $s\succcurlyeq a$) and, if yes, compute the simple
     element $a^{-1}s$ (respectively $s\: a^{-1}$).
 \vspace{-\itemsep}\vspace{-\parsep}
 \end{namedenumerate}
\end{assumption}

We further assume that elements of $G$ are stored as products (sequences) of simple elements or inverses of simple elements.
Then, two elements consisting of at most $k$ such factors can be multiplied at
a cost of $O(k)$ simply by concatenating the corresponding sequences.

We remark that we also could have considered the following additional basic operations:
 \begin{namedenumerate}[Op]\em\vspace{-\topsep}
  \item Given a simple element $s$, test whether $s=1$.
  \item Given two simple elements $s$ and $t$, test whether $s=t$.
  \item Given an atom $a\in\mathcal A$ and a simple element $s$, test whether
   $sa$ (resp.~$as$) is simple and, if yes, compute the simple element $sa$ (resp.~$as$).
 \end{namedenumerate}

However, if $s$ is a simple element, then $s=1$ is equivalent to $a_i\not\preccurlyeq s$
for all $i=1,\dots,\lambda$, where the latter condition can be tested using the operation (Op) at
most $\lambda$ times.  Hence, (Op1) can be realised in terms of (Op) at a cost of $O(C\lambda)$.
We will moreover see below that
(Op2) and (Op3) can be realised in terms of (Op) at a cost of $O(C\lambda\wordlen{\Delta})$.
While doing so may not yield the most efficient ways of realising (Op1), (Op2) and (Op3),
it does not change the complexities of the algorithms we consider.

We remark that the operations (Op) and (Op3) can be realised at equal cost in many Garside groups;
this is the case for braid groups, for instance.  However, as we are working with a generic Garside
group of finite type, we want to keep our assumptions to the minimum.
We moreover mention that one could use (Op3) as basic operation instead of (Op):  if the cost of (Op3) is $O(C)$, then one can test at a cost of $O(C\lambda)$ whether a simple element is equal to $\Delta$ and the operations (Op) and (Op2) can be realised in terms of (Op3) at a cost of $O(C\lambda\wordlen{\Delta})$; the map $\partial$ induces a duality between this situation and the situation from Assumption~\ref{A:BasicOperations}.
Finally, note that (Op1) can be realised in terms of (Op3) at a cost of $O(C\lambda)$, if $\Delta$ is the lcm of the atoms of $G$: in this case, $\partial(s)=\Delta$ is equivalent to $a_i\preccurlyeq \partial(s)$ for all $i=1,\dots,\lambda$, that is,
$s=1$ is equivalent to $sa_i\in[1,\Delta]$ for all $i=1,\dots,\lambda$.

An important remark concerning the algorithms below is the following: One of the most
frequently used operations consists of determining an atom $a$ such that $a\preccurlyeq s$, given
a nontrivial simple element $s$. If the simple elements are stored as products of atoms,
this operation has a cost of~$O(1)$. However, if the simple
elements are stored in a different way, it is possible that the only way to find such an
atom is to check whether $a\preccurlyeq s$ for every $a\in\mathcal A$, until
the answer is positive. This has time complexity $O(C\lambda)$. Therefore, in the
algorithms below we will sometimes write `Take an atom $a\preccurlyeq s$', and we will
assume that this operation has a cost of $O(C\lambda)$, although the reader should notice
that the actual cost could be only $O(1)$ in some situations.

The first computations which we will express in terms of the basic operations are computing left
and right complements of simple elements and conjugation of simple elements by $\Delta$ or $\Delta^{-1}$.  We will also see a generic way of performing the operations (Op2) and (Op3).  The following algorithm underlies all of these:

\begin{minipage}{\textwidth}
\fbox{{\bf Algorithm to compute the right complement of a simple element}}

\begin{longtable}{lp{12.5cm}}
{\bf Input:} & A simple element $s$. \\
{\bf Output:} & The simple element $\partial(s)=s^{-1}\Delta$.
\end{longtable}

\begin{enumerate}
\item Set $d=\Delta$.
\item While $s\neq 1$ do:
  \begin{enumerate}
  \item Take an atom $a\preccurlyeq s$.
  \item Set $d= a^{-1} d$ and $s= a^{-1}s$.
  \end{enumerate}
\item Return $d$.
\end{enumerate}
\end{minipage}

At most $\wordlen{\Delta}$ passes through the loop are required and the costs of the test $s\ne 1$,
step~2\,(a) and step~2\,(b) are $O(C\lambda)$, $O(C\lambda)$ and $O(C)$, respectively. Hence, the
complexity of this algorithm is $O(C\lambda\wordlen{\Delta})$. Notice that
$\partial^{-1}(s)=\Delta \: s^{-1}$ can be computed in the same way,
replacing $\preccurlyeq$ by $\succcurlyeq$ and multiplying with $a^{-1}$ on the right instead of on
the left.
The given algorithm can also be used to compute
$\tau(s)=\partial^2(s)$ or $\tau^{-1}(s)=\partial^{-2}(s)$, so all these operations have
a cost of $O(C\lambda\wordlen{\Delta})$.

Given a simple element $s$ and an atom $a$, one can determine whether $sa$ is simple by computing
$\partial(s)$ with the above algorithm and checking whether $a\preccurlyeq \partial(s)$, where the
latter step has a cost of $O(C)$ by Assumption~\ref{A:BasicOperations}.  Moreover, if $sa$ is
simple, one can compute $sa=\partial^{-1}(a^{-1}\partial(s))$.
That is, we can perform operation (Op3) that way.
Similarly, one can determine whether $as$ is simple by checking whether
$\partial^{-1}(s)\succcurlyeq a$ and, if it is, one can compute
$as=\partial(\partial^{-1}(s) a^{-1})$.  All these operations have a cost of
$O(C\lambda\wordlen{\Delta})$.

Next, we will describe the lattice operations on simple elements, which are important for computing
normal forms of elements.

\begin{minipage}{\textwidth}
\fbox{{\bf Algorithm to compute the greatest common divisor of two simple elements}}

\begin{longtable}{lp{12.5cm}}
{\bf Input:} & Two simple elements $s$ and $t$. \\
{\bf Output:} & The simple element $s\wedge t$.
\end{longtable}

\begin{enumerate}
\item Set $i=1$ and $d=\Delta$.
\item While $i\le\lambda$ do:
  \begin{enumerate}
  \item If $a_i\preccurlyeq s$ and $a_i\preccurlyeq t$, then
  \item \qquad set $d=a_i^{-1}d$, set $s=a_i^{-1}s$, set $t=a_i^{-1}t$ and set $i=1$,
  \item[] else
  \item \qquad set $i=i+1$.
  \end{enumerate}
\item Return $\partial^{-1}(d)$.
\end{enumerate}
\end{minipage}

The tests in step~2\,(a) and the operations in step 2\,(b) have a cost of $O(C)$, step~3 has a
cost of $O(C\lambda\wordlen{\Delta})$, and all remaining operations have a cost
of $O(1)$.  As step~2\,(b) is executed at most $\wordlen{\Delta}$ times, with at most $\lambda$
passes through the while loop between two consecutive executions, the cost of step~2 is
$O(C\lambda\wordlen{\Delta})$, so the complexity of the algorithm is also
$O(C\lambda\wordlen{\Delta})$.
Note that finding the atoms which are common divisors of $s$ and $t$ is critical for the complexity
of the algorithm.  Thus, even if step~3 was avoided by making use of a realisation of (Op3) with a
cost of $O(C)$, the complexity of the algorithm would not improve.

By symmetry, one can similarly compute the greatest common divisor $s \rwedge t$ with respect to
$\succcurlyeq$.

Least common multiples of simple elements with respect to $\preccurlyeq$ or $\succcurlyeq$ can now
be computed using the following formulae, which can easily be seen to hold:
$$
   s\vee t = \partial^{-1}\left(\partial(s) \rwedge \partial(t) \right),
   \hspace{2cm}  s\rvee t =  \partial\left( \partial^{-1}(s) \wedge \partial^{-1}(t) \right).
$$
Therefore, computing $s\vee t$ or $s\rvee t$ also takes time
$O(C\lambda\wordlen{\Delta})$.

As $s=t$ is equivalent to $s=s\wedge t=t$, we can use the following modification of the above algorithm to test whether two simple elements are equal, that is, perform operation (Op2).

\begin{minipage}{\textwidth}
\fbox{{\bf Algorithm to test whether two simple elements are equal}}

\begin{longtable}{lp{12.5cm}}
{\bf Input:} & Two simple elements $s$ and $t$. \\
{\bf Output:} & The truth value of $s=t$.
\end{longtable}

\begin{enumerate}
\item Set $i=1$ and $d=\Delta$.
\item While $i\le\lambda$ do:
  \begin{enumerate}
  \item If $a_i\preccurlyeq s$ and $a_i\preccurlyeq t$, then
  \item \qquad set $d=a_i^{-1}d$, set $s=a_i^{-1}s$, set $t=a_i^{-1}t$ and set $i=1$,
  \item[] else
  \item \qquad set $i=i+1$.
  \end{enumerate}
\item If $s=1$ and $t=1$, then return true, else return false.
\end{enumerate}
\end{minipage}

The cost of step~3 is
$O(C\lambda)$; all other steps are as before.  Hence, the complexity of the algorithm is
$O(C\lambda\wordlen{\Delta}$).
This implies, in particular, that two elements of canonical length at most $k$ whose (left or right) normal forms are known, can be compared at a cost of $O(C\lambda k\wordlen{\Delta})$ by comparing their infima (at a cost of $O(1)$) and at most $k$ pairs of simple elements.

The following algorithm computing the local sliding of a pair of simple elements is also just a small modification of the algorithm as the one computing the gcd of two simple elements:

\begin{minipage}{\textwidth}
\fbox{{\bf Algorithm to compute the local sliding of a pair of simple elements}}

\begin{longtable}{lp{12.5cm}}
{\bf Input:} & Two simple elements $s$ and $t$. \\
{\bf Output:} & The simple elements $s (\partial(s)\wedge t)$ and $(\partial(s)\wedge
t)^{-1}t $.
\end{longtable}

\begin{enumerate}
\item Set $i=1$ and $s'=\partial(s)$.
\item While $i\le\lambda$ do:
  \begin{enumerate}
  \item If $a_i\preccurlyeq s'$ and $a_i\preccurlyeq t$, then
  \item \qquad set $d=a_i^{-1}d$, set $s'=a_i^{-1}s'$, set $t=a_i^{-1}t$ and set $i=1$,
  \item[] else
  \item \qquad set $i=i+1$.
  \end{enumerate}
\item Return $\partial^{-1}(s'),\;t$.
\end{enumerate}
\end{minipage}

The cost of step~1 is
$O(C\lambda\wordlen{\Delta})$; all other steps are identical.  Hence, the local
sliding of a pair of simple elements can also be computed at a cost of
$O(C\lambda\wordlen{\Delta})$.

Knowing how to compute local slidings, one can use the standard algorithms to
compute the left or right normal form of any element (see  \S\ref{SS:Garside}), based on the
following well-known result.

\pagebreak

\begin{proposition}[see, for example, {\cite[Props.~3.1 and 3.3]{Charney}} or \cite{Epstein}]
\label{P:normal_form}
Let $s_1,\dots,s_{k}$ and $s_0',s_{k+1}'$ be simple elements such that the product
$s_1\cdots s_{k}$ is in left normal form as written.
\begin{enumerate}\vspace{-\topsep}

\item  Consider the product $s_0's_1\cdots s_k$.  For $i=1,\ldots,k$ apply a local sliding to the
pair $s_{i-1}' s_i$, that is, let $t_i= \partial(s_{i-1}') \wedge s_i$ and define
$s_{i-1}''=s_{i-1}'t_i$ and $s_i'= t_i^{-1}s_i$. Finally define $s_k''=s_k'$. Then,
$s_0''\cdots s_k''$ is the left normal form of $s_0's_1\cdots s_k$ (where possibly
$s_0''=\Delta$ or $s_k''=1$).

\item  Consider the product $s_1\cdots s_k s_{k+1}'$. For $i=k,\ldots,1$ apply a
local sliding to the pair $s_i s_{i+1}' $, that is, let $t_i= \partial(s_i) \wedge s_{i+1}'$
and define $s_{i}'=s_{i}t_i$ and $s_{i+1}''= t_i^{-1}s_{i+1}'$. Finally define $s_1''=s_1'$. Then,
$s_1''\cdots s_{k+1}''$ is the left normal form of $s_1\cdots s_k s_{k+1}'$ (where possibly
$s_1''=\Delta$ or $s_{k+1}''=1$).

\end{enumerate}
\end{proposition}

Given an element $x$ written as a product of $k$ simple elements or inverses of simple elements,
the left normal form of $x$ can be obtained as follows.  First, one replaces each inverse $s^{-1}$
of a simple element with $\Delta^{-1}\partial^{-1}(s)$;  at most $k$ replacements are necessary and
each replacement has a cost of $O(C\lambda\wordlen{\Delta})$.  Then, one collects all
appearances of $\Delta$ or $\Delta^{-1}$ on the left hand side, applying $\tau$ or $\tau^{-1}$ as
required, so that the element will be written as $\Delta^q s_1\cdots s_k$, where each $s_i$ is a
simple element;  the number of applications of $\tau$ or $\tau^{-1}$ is bounded by $k(k-1)/2$ and
each application has a cost of $O(C\lambda\wordlen{\Delta})$.  Finally, one applies local
slidings to every pair of consecutive simple elements until every pair is left weighted; it follows
from Proposition~\ref{P:normal_form} that at most $k(k-1)/2$ local slidings are required, each at a
cost of $O(C\lambda\wordlen{\Delta})$.  Therefore, the complexity of computing the left
normal form of $x$ is $O(C\lambda k^2\wordlen{\Delta})$. Computing right normal forms
is analogous and has the same complexity.

Note, however, that if the left normal form (resp.~the right normal form) of $x$ is known and $s$ is a simple element, then the left normal forms (resp.~the right normal forms) of $xs$, $sx$, $xs^{-1}$, $s^{-1}x$, $x^s$ and $x^{s^{-1}}$ can be computed at a cost of
$O(C\lambda k\wordlen{\Delta})$, where $k=\ell(x)$:  the number of applications of
$\tau$ or $\tau^{-1}$ is bounded by $k$ and only $O(k)$ local slidings are required by
Proposition~\ref{P:normal_form}.

We now show how to compute the gcd of two arbitrary elements $a$ and $b$, given as products of
simple elements and inverses of simple elements with at most $k$ factors. First, we write
them in left normal form, say $\Delta^{p}a_1\cdots a_r$ and $\Delta^q b_1\cdots b_t$. If
we denote $m=\min\{p,q\}$, we can consider $a'=\Delta^{-m} a$ and $b'=\Delta^{-m}b$.
Notice that $a'$ and $b'$ are positive elements, and one of them has infimum zero.
Since $a\wedge b= \Delta^m a' \wedge \Delta^m b' = \Delta^m (a'\wedge b')$, it is sufficient to
know how to compute gcds of positive elements and we will hence detail the algorithm to
compute $a\wedge b$ assuming $a$ and $b$ are positive; the cost of reducing to this case by
computing the normal forms of $a$ and $b$ is $O(C\lambda k^2\wordlen{\Delta})$.  We
remark that, if the left normal form of a positive element $a$ is known, then $a\wedge \Delta$ is
also known, since it is precisely the first factor in its left normal form (which may be $\Delta$).

\begin{minipage}{\textwidth}
\fbox{{\bf Algorithm to compute the greatest common divisor of two positive elements}}

\begin{longtable}{lp{12.5cm}}
{\bf Input:} & Two positive elements $a$ and $b$. \\
{\bf Output:} & The element $a\wedge b$.
\end{longtable}

\begin{enumerate}

\item Set $u=\Delta$, $a'=a$, $b'=b$ and $d=1$.

\item While $u\neq 1$ do:

\begin{enumerate}

\item Compute the left normal forms of $a'$ and $b'$.

\item Set $s=a'\wedge \Delta$ and $t=b'\wedge \Delta$.

\item Set $u=s\wedge t$.

\item Set $d=du$, set $a'=u^{-1}a'$ and $b'= u^{-1}b'$.

\end{enumerate}

\item Return $d$.

\end{enumerate}
\end{minipage}

\pagebreak

Since $a$ and $b$ are positive, one has $(a\wedge b)\wedge 1 = 1$.  It is then easy to see by
induction that after the $i$-th pass through the while loop one has
$d=(a\wedge b)\wedge \Delta^i$.  Hence, if $a$ and $b$ are given as products of simple elements and
inverses of simple elements with at most $k$ factors, the number of repetitions of the while loop
is bounded by $k+1$.
The cost of step~2\,(a) in the first pass through the while loop is
$O(C\lambda k^2\wordlen{\Delta})$, but in all subsequent passes, the cost is
$O(C\lambda k\wordlen{\Delta})$ by Proposition~\ref{P:normal_form}.
As the costs of steps~2\,(b), 2\,(c) and 2\,(d) are $O(1)$,
$O(C\lambda\wordlen{\Delta})$ and $O(k)$, respectively, the complexity of the
algorithm hence is $O(C\lambda k^2\wordlen{\Delta})$.
Computing the right gcd $a\rwedge b$ is analogous and has the same complexity.

One can now compute the least common multiple of two elements $a$ and $b$, given as products of
simple elements and inverses of simple elements with at most $k$ factors, as follows. Compute the
normal forms of $a$ and $b$ and let $m=\max\{\sup(a),\sup(b)\}$. The elements $a^{-1} \Delta^m$
and $b^{-1}\Delta^m$ are both positive, whence we can compute
$d= (a^{-1}\Delta^m) \rwedge (b^{-1}\Delta^m)$ using (the right version of) the
algorithm above. Then, $a\vee b = (a^{-1}\rwedge b^{-1})^{-1} =
\Delta^m((a^{-1}\Delta^m)\rwedge (b^{-1}\Delta^m))^{-1} = \Delta^{m} d^{-1}$.
This cost of this computation is dominated by computing $d$ as the right gcd of $a^{-1} \Delta^m$
and $b^{-1}\Delta^m$ which has cost $O(C\lambda k^2\wordlen{\Delta})$.  Thus, the
complexity of computing the lcm $a\vee b$ is $O(C\lambda k^2\wordlen{\Delta})$.
Computing the right lcm $a\rvee b$ is analogous and has the same complexity.

The computations of the preferred prefix and the cyclic sliding of an element can now be
done just by applying the definitions, since we already know how to perform all
operations that occur. For instance, in order to compute the preferred prefix of an
element $x$, given as a product of simple elements and inverses of simple elements with $k$ factors,
one first computes the left normal form of $x=\Delta^p x_1\cdots x_r$, which takes time
$O(C\lambda k^2\wordlen{\Delta})$. Then one applies the formula given in
Definition~\ref{D:preferred prefix}, namely $\p(x)=\iota(x)\wedge \partial(\varphi(x))$. Since
$\iota(x)=\tau^{-p}(x_1)$ with $|p|\le k$ and $\varphi(x)=x_r$, the complexity
of computing $\p(x)$ from the normal form of $x$ is $O(C\lambda k\wordlen{\Delta})$.
The normal form of $\s(x)=x^{\p(x)}$ can then be computed in
$O(C\lambda k\wordlen{\Delta})$.  Thus, the cost of applying a cyclic sliding is
dominated by the cost of computing the normal form, that is, applying a cyclic sliding  has
complexity $O(C\lambda k^2\wordlen{\Delta})$.  Note that if the normal form of $x$ is
known, then $\p(x)$ and the normal form of $\s(x)$ can be obtained at a cost of
$O(C\lambda k\wordlen{\Delta})$.

The transport of an element $\alpha$ at an element $x$ is given by the formula
$\alpha^{(1)}=\p(x)^{-1} \alpha \p(x^\alpha)$. If $x$ and $\alpha$ are given as products of simple
elements and inverses of simple elements with at most $k$ factors, then $\alpha^{(1)}$ can be
computed with the above formula in time $O(C\lambda k^2\wordlen{\Delta})$ by the
arguments from the previous paragraph. In other words, applying a transport has the same complexity
as computing a normal form.
Note that if the normal form of $x$ is known and $\alpha$ is simple, then the normal form of
$x^\alpha$ can be obtained at a cost of $O(C\lambda k\wordlen{\Delta})$, whence
$\alpha^{(1)}$ can be computed at a cost of $O(C\lambda k\wordlen{\Delta})$ by the
arguments above.

Computing the preferred suffix, applying a cyclic right sliding and applying right transport are
analogous and the complexities are the same as for the left versions discussed above.

Finally, the pullback of a positive element $s$ at an element $y$, with the hypotheses and the
notation of Proposition~\ref{P:pullback_computation}, is
$s_{(1)}= \left( \p(z) \ s \ \rp(y^s)^{-1} \right) \;\vee\; 1$;  we assume that we also know the
element~$z$. If $y$, $z$ and $s$ are given as products of simple elements and inverses of simple
elements with at most $k$ factors, then $s_{(1)}$ can be computed in time
$O(C\lambda k^2\wordlen{\Delta})$ using the operations described above.
If $s$ is simple and if the left normal form of $z$ and the right normal form of $y$ are known, then
$\p(z) \ s \ \rp(y^s)$ can be computed at a cost of
$O(C\lambda k\wordlen{\Delta})$ and, since this product involves only 3 simple
factors, the subsequent computation of the lcm has a cost of
$O(C\lambda \wordlen{\Delta})$, whence in this case $s_{(1)}$ can be obtained at a
cost of $O(C\lambda k\wordlen{\Delta})$.
Computing the right pullback $s_{\rt{(1)}}$ is analogous and has the same complexity.

\pagebreak

Summarising the results obtained in this section, we have:
 \begin{theorem}\label{T:BasicOperations}
  Let $G$ be a Garside group of finite type with Garside element $\Delta$ and set of atoms
  $\mathcal A=\{a_1,\ldots,a_{\lambda}\}$ for which Assumption~\ref{A:BasicOperations} is
  satisfied.  Moreover, let $a$ be an atom of $G$, let $s$ and $t$ be simple elements of $G$ and
  let $x$, $y$ and $\alpha$ be elements of $G$, given as products of simple elements or inverses of
  simple elements with at most $k$ factors.
  \begin{enumerate}\vspace{-\topsep}
   \item The following operation can be performed in $O(C\lambda)$:
   \begin{itemize}
    \item Test whether $s=1$.
   \end{itemize}
   \item The following operations can be performed in $O(C\lambda\wordlen{\Delta})$:
   \begin{itemize}
    \item Test whether $s=t$.
    \item Compute $\partial(s)$, $\partial^{-1}(s)$, $\tau(s)$ or $\tau^{-1}(s)$.
    \item Test whether the product $as$ is simple and, if so, compute $as$.
    \item Test whether the product $sa$ is simple and, if so, compute $sa$.
    \item Compute $s\wedge t$, $s\rwedge t$, $s\vee t$ or $s\rvee t$.
    \item Perform a local (left or right) sliding on the product $s\cdot t$.
   \end{itemize}
   \item The following operations can be performed in $O(C\lambda k\wordlen{\Delta})$:
   \begin{itemize}
    \item Test whether $x=y$, if the left normal forms or the right normal forms of $x$ and $y$ are
          known.
    \item Compute the left normal form \textup{[}resp.~the right normal form\textup{]} of $xs$,
          $sx$, $xs^{-1}$, $s^{-1}x$, $x^s$ or $x^{s^{-1}}$, if the left normal form
          \textup{[}resp.~the right normal form\textup{]} of $x$ is known.
    \item Compute $\p(x)$ or $\s(x)$ \textup{[}resp.~$\rp(x)$ or $\rs(x)$\textup{]}, if the left
          normal form \textup{[}resp.~the right normal form\textup{]} of $x$ is known.
    \item Compute the left transport $s^{(1)}$ \textup{[}resp.~the right transport
          $s^{\rt{(1)}}$\textup{]} of $s$ at $x$, if the left normal form \textup{[}resp.~the right
          normal form\textup{]} of $x$ is known.
    \item Compute the left pullback $s_{(1)}$ \textup{[}resp.~the right pullback
          $s_{\rt{(1)}}$\textup{]} of $s$ at $x$, if it is defined and if the right normal form
          \textup{[}resp.~the left normal form\textup{]} of $x$ and the left normal form
          of the element $z\in\SC(x)$ satisfying $\s(z)=x$ \textup{[}resp.~the right normal form of
          the element $z\in\SC^{\Lsh}(x)$ satisfying $\rs(z)=x$\textup{]} are known.
   \end{itemize}
   \item The following operations can be performed in $O(C\lambda k^2\wordlen{\Delta})$:
   \begin{itemize}
    \item Compute the left normal form of $x$ or the right normal form of $x$.
    \item Compute $x\wedge y$, $x\rwedge y$, $x\vee y$ or $x\rvee y$.
    \item Compute $\p(x)$, $\rp(x)$, $\s(x)$ or $\rs(x)$.
    \item Compute the left transport $\alpha^{(1)}$ of $\alpha$ at $x$ or the right transport
          $\alpha^{\rt{(1)}}$ of $\alpha$ at $x$.
    \item Compute the left pullback $\alpha_{(1)}$ \textup{[}resp.~the right pullback
          $\alpha_{\rt{(1)}}$\textup{]} of $\alpha$ at $x$, if it is defined and if the element
          $z\in\SC(x)$ satisfying $\s(z)=x$ \textup{[}resp.~the element $z\in\SC^{\Lsh}(x)$
          satisfying $\rs(z)=x$\textup{]} is known.
   \end{itemize}
  \end{enumerate}
 \end{theorem}

\subsection{Complexity of the new algorithms}\label{S:ComplexityNew}

Knowing the computational cost of the basic operations, we can now analyse the complexity of the
algorithms for computing $\SC(x)$ from Section~\ref{SS:algorithm}.  Firstly, we define some bounds
which will be used in the sequel.

\begin{notation}\label{N:bounds}
Let $x$ be an element of $G$ given as a product of simple elements or inverses of simple elements
with at most $k$ factors.\vspace{-\topsep}
 \begin{description}
  \item[{[Distance to cyclic sliding repetition]}] Let $T$ be an integer such that there exist
        two integers $0\le i < j \le T$ satisfying $\s^i(x)=\s^j(x)$.
  \item[{[Length of sliding circuits]}] Let $M$ be an integer such that for any element
        $z\in\SC(x)$ there exists a positive integer $N\le M$ with $\s^N(z)=z$.
  \item[{[Distance to transport repetition]}] Let $R$ be an integer such that for any
        element $z\in\SC(x)$ and any simple element $s$ satisfying $z^s\in\SSS(x)$ there exist two
        integers $0\le i < j \le R$ satisfying $s^{(iN)}=s^{(jN)}$, where
        $\s^N(z)=z$ and $s^{(m)}$ denotes $m$-fold transport at $z$ for $m\in\NN$.
 \end{description}
\end{notation}

\begin{remark}
\label{R:obvious_bounds}\em
It is easy to see that integers $T$, $M$ and $R$ as above exist and to give some obvious (but very
crude) upper bounds for them:  By Corollary~\ref{C:sliding_reaches_period}, iterated cyclic sliding
becomes periodic, so $T$ as above exists.  Indeed, it follows from
Proposition~\ref{P:sliding_reaches_SSS} that $\s^{m}(x)\in\SSS(x)$ for all
$m\ge k\wordlen{\Delta}$.  As $|\SSS(x)|\le |[1,\Delta]|^k$, one can choose
$T\le k\wordlen{\Delta} + \;|[1,\Delta]|^k$.
Moreover, as $\SC(x)\subseteq\SSS(x)$ is finite, $M$ as above exists and one can choose
$M\le|\SC(x)|$.  (Hence, in particular, $M\le|[1,\Delta]|^k$.)  Finally, by
Proposition~\ref{P:transport}~(5), transports of simple elements are simple.  Since $G$ is of
finite type, $R$ as above exists and one can choose $R\le|[1,\Delta]|$.
\end{remark}

\begin{lemma}\label{L:pullback_repetition}
Let $x\in G$, $z\in\SC(x)$, and let $s$ be a simple element such that $z^s\in\SSS(x)$.
If $N$, $i$, $j$ and $K$ are integers such that $\s^N(z)=z$, $0\le i < j \le K$ and
$\left(s_{(KN)}\right)^{(iN)} = \left(s_{(KN)}\right)^{(jN)}$, where $t^{(m)}$ denotes $m$-fold
transport of $t$ at $z$ and $t_{(m)}$ denotes $m$-fold pullback of $t$ at $z$ for $m\in\NN$,
then $s_{(KN)}=s_{((K+j-i)N)}$.
\end{lemma}

\begin{proof}
By Lemma~\ref{L:iterated_pullback} we have
$s_{(KN-iN)} \preccurlyeq \left(\left(s_{(KN-iN)}\right)_{(iN)}\right)^{(iN)}
 = \left(s_{(KN)}\right)^{(iN)} = \left(s_{(KN)}\right)^{(jN)}$.
Again using Lemma~\ref{L:iterated_pullback}, we obtain
$\left(s_{(KN-iN)}\right)_{(jN)} \preccurlyeq s_{(KN)}$, that is,
$s_{((K+j-i)N)}\preccurlyeq s_{(KN)}$.

Similarly, we have $s_{(KN-jN)} \preccurlyeq \left(\left(s_{(KN-jN)}\right)_{(jN)}\right)^{(jN)}
 = \left(s_{(KN)}\right)^{(jN)} = \left(s_{(KN)}\right)^{(iN)}$ and from this obtain
$\left(s_{(KN-jN)}\right)_{(iN)} \preccurlyeq s_{(KN)}$, that is,
$s_{((K+i-j)N)}\preccurlyeq s_{(KN)}$.  Applying $(j-i)N$-fold pullback to the last statement
yields $s_{(KN)}\preccurlyeq s_{((K+j-i)N)}$ using Lemma~\ref{L:pullback_monotonic}.

Hence, $s_{(KN)}=s_{((K+j-i)N)}$ as we wanted to show.
\end{proof}

\begin{corollary}
\label{C:pullback_repetition}
Consider for $x\in G$ the bounds from Notation~\ref{N:bounds}.
For any element $z\in\SC(x)$ and any simple element $s$ satisfying $z^s\in\SSS(x)$ there exist two
integers $0\le i < j \le 2R$ satisfying $s_{(iN)}=s_{(jN)}$, where $\s^N(z)=z$ and $s_{(m)}$
denotes $m$-fold pullback at $z$ for $m\in\NN$.
\end{corollary}

\begin{proof}
By the choice of $R$ there are integers $0\le i'<j'\le R$ such that
$\left(s_{(RN)}\right)^{(i'N)} = \left(s_{(RN)}\right)^{(j'N)}$. We then have
$s_{(RN)}=s_{((R+j'-i')N)}$ by Lemma~\ref{L:pullback_repetition}.  Setting $i=R$ and $j=R+j'-i'$,
we have $0\le i<j\le 2R$ and $s_{(iN)}=s_{(jN)}$ as desired.
\end{proof}

\begin{proposition}\label{P:complexity_algo_1}
Let $G$ be a Garside group of finite type with Garside element $\Delta$ and $\lambda$ atoms, and
let $x$ be an element of $G$ given as a product of simple elements or inverses of simple elements
with at most $k$ factors.  Using the bounds from Notation~\ref{N:bounds}, the complexity of
Algorithm 1 is $O(C\lambda k(k+T)\wordlen{\Delta})$.
\end{proposition}

\begin{proof}
Observe that $\ell(\s^i(x))\le k$ for all non-negative integers $i$.
In particular, the normal forms of two such elements can be compared at a cost of
$O(C\lambda k\wordlen{\Delta})$ by Theorem~\ref{T:BasicOperations}.
Note further that a hash function depending on all factors in the normal
form can be computed at a cost of $O(Ck)$, if the normal form is known.  We use a sufficiently large
hash table, together with this hash function, to store the trajectory $\mathcal T$ in step~2.  If
the normal form of an element $y$ with $\ell(y)\le k$ is known, testing whether $y\in\mathcal T$
(and storing it if it is not) then has a cost of $O(C\lambda k\wordlen{\Delta})$.

We initially compute the normal form of $x$ at a cost of $O(C\lambda k^2\wordlen{\Delta})$.
Step~1 has a cost of $O(1)$.  Step~3 and each pass through the while loops in step~2 and step~4 have
a cost of $O(C\lambda k\wordlen{\Delta})$ by Theorem~\ref{T:BasicOperations}.  The
number of passes through the while loops is bounded by T. Step~5 has a cost of $O(T)$.  Hence the
claim holds.
\end{proof}

\begin{proposition}\label{P:complexity_algo_2}
Let $G$ be a Garside group of finite type with Garside element $\Delta$ and $\lambda$ atoms, and
let $v$ be an element of $G$ given as a product of simple elements or inverses of simple elements
with at most $k$ factors.  If the left and right normal forms of $v$ are known, then, using the
bounds from Notation~\ref{N:bounds}, the complexity of Algorithm 2 is
$O\big(C\lambda^2 k\wordlen{\Delta}(\wordlen{\Delta}+RM)\big)$.
\end{proposition}

\begin{proof}
In step~1, we perform $N\le M$ times the following operations: apply a cyclic sliding to an element
whose left and right normal forms are known, compute the left normal form and the right normal form
of the result and compare it to $v$; each of these has a cost of
$O(C\lambda k\wordlen{\Delta})$ by Theorem~\ref{T:BasicOperations}. Hence, the cost
of step~1 is $O(C\lambda kM\wordlen{\Delta})$.

Step~2 has a cost of $O(\lambda)$; we store the set $\mathcal A_v$ as a list and the set
$Atoms\subseteq\{a_1,\dots,a_\lambda\}$ as an array of $\lambda$ flags.

Steps~3\,(a) to (e) are executed $\lambda$ times.  Step~3\,(a) has a cost of $O(1)$; the costs of
the remaining steps are as follows:

For step~3\,(b) note that at any time we have
$a_t\preccurlyeq s\preccurlyeq\rho_{a_t}\preccurlyeq\Delta$, so $s$ is simple.  In particular,
$\sup(v^s)-\sup(v)\in\{0,1\}$.  As the right normal form of $v$ is known and $s$ is simple, the
right normal form of $v^s$ can be computed at a cost of
$O(C\lambda k\wordlen{\Delta})$ by Theorem~\ref{T:BasicOperations}.
By Proposition~\ref{P:alg_rho_s_lcm}, we can obtain the element
$1\vee v^s\Delta^{-\sup(v)}$ from the right normal form of $v^s$ at a cost of $O(1)$: it is the
leftmost factor in the right normal form if $\sup(v^s)=\sup(v)+1$, and it is trivial if
$\sup(v^s)=\sup(v)$.  In the same way, we can obtain $1\vee (v^s)^{-1}\Delta^{\inf(v)}$ from the
right normal form of $(v^s)^{-1}$.  Observe that the right normal form of $(v^s)^{-1}$ is related
to the right normal form of~$v^s$: the leftmost factor in the right normal form of $(v^s)^{-1}$ can
be obtained from the rightmost non-$\Delta$ factor in the right normal form of $v^s$ by applying
the map $\partial$ or $\partial^{-1}$ at most $2k+1$ times, that is, at a cost of
$O(C\lambda k\wordlen{\Delta})$ by Theorem~\ref{T:BasicOperations}.  As both
$1\vee v^s\Delta^{-\sup(v)}$ and $1\vee (v^s)^{-1}\Delta^{\inf(v)}$ are simple, so is their lcm. 
In particular, computing the lcm and the final multiplication (which is a local sliding) each have
a cost of $O(C\lambda \wordlen{\Delta})$ by Theorem~\ref{T:BasicOperations}.
Hence, since the number of passes through the while loop is at most $\wordlen{\Delta}$ by
Proposition~\ref{P:alg_rho_s}, step~3\,(b) has a cost of
$O(C\lambda k\wordlen{\Delta}^2)$.

In step~3\,(c) the initial test $a_t\preccurlyeq\p(v)$ has a cost of $O(C)$.  We can store the
simple elements $s_{(iN)}$ ($i=1,2\dots$) in a sufficiently large hash table, using the hash
function from Assumption~\ref{A:BasicOperations}.  Testing whether $s_{(iN)}$ has already
occurred (and storing it if not) then has a cost of $O(C\lambda\wordlen{\Delta})$.  Since the left and right normal
forms of all elements in the sliding circuit of $v$ are known from step~1, each pullback can be
computed at a cost of $O(C\lambda k\wordlen{\Delta})$ by
Theorem~\ref{T:BasicOperations}.  As the number of pullbacks which need to be computed is bounded
by $2RM$ by Corollary~\ref{C:pullback_repetition}, the cost of step~3\,(c) hence is
$O(C\lambda kRM\wordlen{\Delta})$.

By the same arguments, step~3\,(d) has a cost of $O(C\lambda kRM\wordlen{\Delta})$,
since each transport can be computed at a cost of $O(C\lambda k\wordlen{\Delta})$ and
the number of transports which need to be computed is bounded by $RM$.

The test in the outer if statement in step~3\,(e) has a cost of $O(CR)$, whereas the test in the
if statement in step~3\,(e)\,i has a cost of $O(C\lambda)$, since testing whether $a_k\in Atoms$ has a cost of~$O(1)$.  As the remaining operations in step~3\,(e)\,i have a cost of $O(1)$, the cost of step~3\,(e) is $O(C(\lambda+R))$.

Hence, the complexity of Algorithm~2 is
$O(C\lambda^2 k\wordlen{\Delta}(\wordlen{\Delta}+RM))$ as claimed.
\end{proof}

\begin{theorem}\label{T:complexity_algo_3}
Let $G$ be a Garside group of finite type with Garside element $\Delta$ and $\lambda$ atoms, and
let $x$ and $y$ be elements of $G$ given as products of simple elements or inverses of simple
elements with at most $k$ factors.  Let $T$, $M$ and $R$ be the maxima of the bounds from
Notation~\ref{N:bounds} for $x$ and $y$, respectively.

The complexity of Algorithm 3 is $O\Big( C\lambda k\wordlen{\Delta}
     \cdot \big( k+T+|\SC(x)|\lambda(\wordlen{\Delta}+RM) \big) \Big)$.
\end{theorem}

\begin{proof}
Observe that $\ell(z)\le k$ for all $z\in\SC(x)$.
In particular, the (left) normal forms of two such elements can be compared at a cost of
$O(C\lambda k\wordlen{\Delta})$ by Theorem~\ref{T:BasicOperations}.
Note further that a hash function depending on all factors in the normal form can be computed at a cost of $O(Ck)$, if the normal form is known.
We use a sufficiently large hash table, together with this hash function, to store the set $\mathcal V$.
More precisely, whenever a new element $v^s\in\SC(x)$ is found, where $s$ is an indecomposable conjugator and $v\in\mathcal V$, we store the left normal form and the right normal form of~$v^s$, as well as the indecomposable conjugator $s$, in the hash table entry for $v^s$.
If the left normal form and the right normal form of $v^s$ are known, testing whether $v^s\in\mathcal V$, and storing all required data if it is not, has a cost of
$O(C\lambda k\wordlen{\Delta})$.
The set $\mathcal V'$ is stored as a list (storing hash table indices instead of actual elements), whence storing or retrieving an element of $\mathcal V'$ has a cost of $O(1)$.
Observe that the conjugating elements $c_v$ for $v\in\mathcal V$ are implicit in the spanning tree structure for $\SCG(x)$ with root $\widetilde{x}$ which is computed: for any $v\in\mathcal V$, the conjugating element~$c_v$ can be obtained by tracing back the path to the root which is given by the indecomposable conjugators stored for every entry in the hash table.  In particular, there is no actual computation of  $c_{v^s}=c_v\cdot s$ in step~3\,(c)\,ii; at most $c_{\widetilde{y}}$ is ever explicitly computed (in step~3\,(c)\,i).

Step~1 has a cost of $O(C\lambda k(k+T)\wordlen{\Delta})$ by Proposition~\ref{P:complexity_algo_1}; this includes computing the left and right normal forms of $\widetilde{x}$ and $\widetilde{y}$.
For step~3\,(c) note that, since the left normal form and the right normal form of $v$ are known, the left normal form and the right normal form of each conjugate~$v^s$ can be computed at a cost of
$O(C\lambda k\wordlen{\Delta})$ by Theorem~\ref{T:BasicOperations}.
Steps~3\,(a), 3\,(d) and 4 have a cost of~$O(1)$.  Steps~2, 3\,(c)\,ii, as well as the test of the
condition in step~3\,(c)\,i have a cost of $O(C\lambda k\wordlen{\Delta})$.
Step~3\,(b) has a cost of $O(C\lambda^2 k\wordlen{\Delta}(\wordlen{\Delta}+RM))$ by
Proposition~\ref{P:complexity_algo_2}.
The body of the while loop in step~3 is executed $|\SC(x)|$ times and the body of the for loop in step~3\,(c) is executed at most $\lambda$ times.
The actual computation of the conjugating element
$c_1\cdot c_{\widetilde{y}}\cdot c_2^{-1}$ in step~3\,(c)\,ii has a cost of $O(T+|\SC(x)|)$, but is
executed at most once.

Thus, the complexity of Algorithm~3 is
\begin{eqnarray*}
  O\Big( C\lambda k(k+T)\wordlen{\Delta} \Big)
  &+& O\Big( |\SC(x)|\cdot C\lambda^2 k\wordlen{\Delta}(\wordlen{\Delta}+RM) \Big) \\
  &+& O\Big( |\SC(x)|\lambda\cdot C\lambda k\wordlen{\Delta} \Big)
  \;\;+\;\; O\Big( T+|\SC(x)| \Big) \\
  \lefteqn{=\;\;
  O\Big( C\lambda k\wordlen{\Delta}
     \cdot \big( k+T+|\SC(x)|\lambda(\wordlen{\Delta}+RM) \big) \Big)
  }
\end{eqnarray*}
as claimed.
\end{proof}

\begin{remark}\em
\label{R:better_bounds}
Unfortunately, the obvious bounds for $T$ and $M$ given in Remark~\ref{R:obvious_bounds} are exponential in $k$.  For the Artin braid groups $B_n$ one has $|[1,\Delta]|=n!$, that is, the above bounds are also exponential in $n$ (or $\wordlen{\Delta}$) for this sequence of Garside groups, as is the bound for $R$ given in Remark~\ref{R:obvious_bounds}.
Moreover, no bound for $|\SC(x)|$ is currently known which is better than the obvious bound $|\SC(x)| \le |\SSS(x)| \le |[1,\Delta]|^k$ (cf.~Remark~\ref{R:obvious_bounds}); the latter again is exponential.
None of these bounds adequately describes the behaviour observed in computer experiments.

We conjecture that there are bounds for $T$, $M$ and $R$ which are polynomial in $k$ and $\wordlen{\Delta}$.  If the elements of $\SC(x)$ are rigid, then one can choose $R=\wordlen{\Delta}$ by \cite[Proposition~4.3 and Corollary~4.4]{GG1}, and obviously $M=1$.  However, even in this case, no realistic bound for $T$ is  known.

The situation for $|\SC(x)|$ is more complicated.  It is shown in \cite{BGG3} that $|\USS(x)|$ grows exponentially in $n$ for periodic elements of the Artin braid groups $B_n$. By \cite[Proposition~5.1]{GG1}, the same is true for $|\SC(x)|$.  Hence, a bound for $|\SC(x)|$ which is polynomial in $k$ and $\wordlen{\Delta}$ cannot be expected in general.  However, it may be possible to establish such a bound for certain classes of elements, for example rigid elements.  For the situation of Artin braid groups, an attempt to reduce the general case to the special case of rigid elements is sketched in \cite{BGG1}.

The problem of finding bounds for $T$, $M$ and $R$ which are polynomial in $k$ and $\wordlen{\Delta}$ and the problem of understanding $|\SC(x)|$ correspond to open problems formulated in \cite{BGG1} in the context of Artin braid groups for ultra summit sets and the cycling and decycling operations.
\end{remark}

\vspace{.3cm}
\noindent {\footnotesize
\begin{minipage}[t]{5.2cm}
{\bf Volker Gebhardt:} \\
School of Computing and Mathematics\\
University of Western Sydney \\
Locked Bag 1797 \\
Penrith South DC NSW 1797, Australia
\\ E-mail: v.gebhardt@uws.edu.au
\end{minipage}
\hfill
\begin{minipage}[t]{5.4cm}
{\bf Juan Gonz\'alez-Meneses:} \\
Dept.~\'{A}lgebra.  Facultad de Matem\'{a}ticas\\
Universidad de Sevilla \\
Apdo.~1160 \\
41080 Sevilla (SPAIN)
\\ E-mail:  meneses@us.es
\\ URL: www.personal.us.es/meneses
\end{minipage}
}

\end{document}